%% file: c.tex
\documentclass[12pt]{amsart}

\usepackage{xcolor}

\usepackage{xcolor}

\newcommand{\MNLxi}{\MM\NN(E(L),\xi)}
\newcommand{\MNLo}{\MM\NN(E(L),\O)}
\newcommand{\hmz}{H^1(M,\zz)}
\newcommand{\hmr}{H^1(M,\rr)}

\newcommand{\hmq}{H^1(M,\qq)}

\newcommand{\MNxi}{\MM\NN(M,\xi)}
\newcommand{\Pxi}{{\wh \PPPP}_\xi}

\input bookman.sty
\boldmath

\usepackage{helvet}

\usepackage{graphics}

\usepackage{graphicx}

\usepackage{tikz}
\usepackage{subcaption}
\usepackage{graphicx}

\usepackage{amssymb}
\usepackage{amsxtra}
\usepackage{amsmath}
\usepackage{mathrsfs}

\usepackage[arrow, matrix]{xy}
\xyoption{frame}

\theoremstyle{plain}

\newtheorem{theo}{Theorem}[section]
\newtheorem{lemm}[theo]{Lemma}
\newtheorem{prop}[theo]{Proposition}
\newtheorem{coro}[theo]{Corollary}
\newtheorem{conj}[theo]{Conjecture}

\theoremstyle{definition}

\newtheorem{defi}[theo]{Definition}
\newtheorem{rema}[theo]{Remark}

\newtheorem{examples}[theo]{Examples}

\newfont{\rmm}{cmr10 scaled 1000}

\newfont{\itt}{cmsl10 scaled 1000}

\newfont{\rM}{cmr10 scaled 1700}

\newcounter{lemma}[section]

\newcounter{tempcounter}

\newcommand{\lb}{\label}

\newcommand{\rrf}[1]{(\ref{#1})}

\input auxxx.tex

\newcommand{\ox}{{\overline {X}} }
\newcommand{\oox}{{{\overline{X}_0}} }
\newcommand{\HN}[2]{ H_*({#1}, {#2})}
\newcommand{\cov}[3]{{#1} \arrr{#2} {#3}}
\newcommand{\arrh}[3]
{
\xymatrix{
{#1} \ar[r]^<<<<{#2}  &{#3}
}
}

\newcommand{\arlh}[3]
{
\xymatrix{
{#1} & \ar[l]_<<<<{#2}  {#3}
}
}

\newcommand{\arrr}[1]
{\arrh {}{#1}{}}

\newcommand{\arrl}[1]
{\arlh {}{#1}{}}

\newcommand{\arl}
{\arrl {}}

\newcommand{\arrto}
{\xymatrix{{} \ar@{|-{>}}[r]  & {} } }

\newcommand{\arrinto}
{\xymatrix{{} \ar@{^{(}->}[r]  & {} } }

\newcommand{\cco}{ chain complex}

\newcommand{\dfox}[2]{\frac {\pr {#1}}{\pr {#2}}}
\newcommand{\tbtc}{2-bridge 2-component link}

\newcommand{\ruru}{\frac {\pr r}{\pr u}}
\newcommand{\rvrv}{\frac {\pr r}{\pr v}}

\newcommand{\ukk}{1+\ldots + u^{k-1}}

\newcommand{\mntt}{Morse-Novikov theory}
\newcommand{\mnn}{Morse-Novikov number}

\begin{document}

\title
[Morse-Novikov  theory for  links]
{Morse-Novikov  theory for  links
}
\author{L. Chen, H. Endo,  A.  Pajitnov}
\address{Department of Mathematics
Tokyo Institute of Technology
2-12-1 Ookayama, Meguro-ku
Tokyo 152-8551
Japan}
\email{chen.l.aj@m.titech.ac.jp}
\address{Department of Mathematics
Tokyo Institute of Technology
2-12-1 Ookayama, Meguro-ku
Tokyo 152-8551
Japan}
\email{endo@math.titech.ac.jp}
\address{Laboratoire Math\'ematiques Jean Leray
UMR 6629,
 Universit\'e de Nantes,
Facult\'e des Sciences,
2, rue de la Houssini\`ere,
44072, Nantes, Cedex}
\email{andrei.pajitnov@univ-nantes.fr}

\newcommand{\mnlx}{\MM\NN(L,\xi)}

\begin{abstract}
For a compact 3-manifold W. Thurston introduced a norm on
the first cohomology group of the manifold.
The unit ball $B$ of this norm is a polyhedron and the
set of cohomology classes that are representable
by fibrations over a circle
is a union of cones on some of the open faces of $B$.
In the present paper we study the fibred faces of the Thurston polyhedra
of  exteriors of links in $S^3$.
Our approach is based on  the non-abelian Novikov homology
associated with the universal covering of the exterior of the link.
We prove in particular that for a \tbtc~  $L$ a cohomology
class $\xi\in H^1(E(L))$ can be represented by a fibration over a circle
\ifff~ its 2-variable Alexander polynomial is $\xi$-monic.
We compute the Morse-Novikov numbers for a majority of
2-component prime links with number of crossings $\leq 8$.
\end{abstract}

\maketitle
\section{Introduction}
\label{s:intro}

In \cite{Thurston}
W. Thurston introduced a norm on the cohomology $H^1(M)$
of any compact 3-manifold $M$. This construction
has many applications in geometry and topology, in particular to
study of fibrations of 3-manifolds over a circle.
Namely, the unit ball $B$ of Thurston's  norm is a polyhedron
and the set of deRham cohomology classes
$\xi\in
H^1(M,\rr)$
representable by
a non-singular 1-form $\o$
%a fibration over $S^1$
is the union of cones on some of
the open faces of $B$.

%In his article  W. Thurston studied these open faces
%for some families of examples.

In the present paper we study the fibred faces of Thurston polyhedron
for complements of links in $S^3$.
We use the fibering obstructions
provided by Morse-Novikov theory.

Let $M$ be a compact manifold, and $\xi\in H^1(M)\approx [M, S^1]$
a cohomology class. Novikov theory provides obstructions for $\xi$ to
be represented by a $\smo$ fibration over a circle (such class will
be called {\it fibered}). Namely, to any regular covering
$p:\bar M \to M$ with structure group $G$ such that $p^*(\xi)$ is
homologous to $0$, one associates the homology module
$
H_*(M,\xi) =
H_*(\Lxi\tens{\L} C_*(\bar M))
$
where $\L=\zz G$ and $\Lxi$ is the Novikov completion of $\L$ \wrt~ $\xi$.
Vanishing of this homology is a necessary condition for $\xi$ to be fibered.
%In the case when $G$ is free abelian this homology is a
%module over a Novikov completetion of Laurent polynomial ring.
%For the case when $M$ is the exterior of an $n$-component link $L$ %in $S^3$
%the corresponding obstructions can be read from the $n$-variable
%Alexander polynomial.

In general vanishing of these obstructions is not sufficient
for fibering. However in dimension $3$,
a profound theorem due to R. Bieri, W. Neumann and R. Strebel
\cite{BNS} and J.-Cl. Sikorav \cite{SikoravThese}
says that  the vanishing of the Novikov homology
associated to the universal covering $\wi M\to M$
is sufficient for fibering.
In the present paper we apply this result to the case when $M$ is the exterior $E(L)$
of a link $L$  in $S^3$.
In Section 1  we deal with the case of \tbtc s.  We prove  that
$\xi$ is fibered \ifff~
the 2-variable Alexander polynomial of
the link in $\xi$-monic. Therefore the fibredness
of $\xi$ is determined in this case by abelian invariants.

In Section 2 we consider connected sums of links.
Let $L=L_1\krest L_2$,
and $\xi\in H^1(E(L))$.
Under a mild restriction on $\xi$ we prove that it is fibered
\ifff~ both its restrictions to the complements
$E(L_1)$
and
$E(L_2)$
are fibered.
This generalizes
a well-known result: the connected sum of knots is fibered \ifff~ each of the knots is fibered
(see D. Gabai's paper \cite{Gabai}). Our proof is a short computation with the Novikov homology, and does not use  foliation theory.

We determine all the fibered cohomology classes
for  prime links with 2 components and at most 8
crossings
(except the link L8a7 for which we have a partial result).
See the tables in Section \ref{s:6-7-8}.

In the second part of the paper we include the question of fibering
of links into more general setting of Morse-Novikov theory of maps to the circle.
Let $M$ be a manifold.
A cohomology class
$\xi\in H^1(M)$ is called {\it regular } if
its restriction to $\pr M$ can be represented by a fibration
over $S^1$.
For a manifold $M$ and a  regular cohomology class
$\xi$
the Morse-Novikov number
$\MM\NN(M,\xi)$
is defined as the minimal number of critical points
of a map $f:M\to S^1$ \sut~ $[f]=\xi$ and $f~|~\pr M$ is a fibration.
Novikov homology provides lower bounds for \mnn s.
For a \cco~ $C_*$ over a ring $R$ denote by $\mu_R(C_*)$ the minimal possible
number of generators of a free \cco~ $D_*$ chain equivalent to $C_*$.
Then we have
\bq\label{f:mn-eq1}
\mu_{\Lxi}\left(\Lxi\tens \L {C_*(\bar M) } \right)\leq \MM\NN(M,\xi).
\end{equation}

These algebraic invariants are easy to handle when $G\approx \zz$ since in this case the ring $\Lxi$ is a principal ideal domain, and
the inequalities
\rrf{f:mn-eq1}
take the familiar Morse-theoretic  form:

\bq\label{f:mn-eq2}
\sum_{i} \left( \wh b_i (M,\xi) +\wh
 q_i (M,\xi) +
  \wh q_{i-1} (M,\xi) \right)
  \leq
\MNxi
\end{equation}
where $\wh b_i, \wh q_i$ denote the Betti and torsion numbers
of the $\Lxi$-module
$\Lxi\tens{\L} H_i(\bar M)$.
In the case when $G\approx \zz^n$ the inequalities
\rrf{f:mn-eq1} can not be reduced to the form
\rrf{f:mn-eq1} but are still moderately accessible.
In the non-abelian case the algebraic invariant
from \rrf {f:mn-eq1} is difficult to compute.

 Observe that  any orientation of $L$ determines
the {\it orientation cohomology class $\O$}, that
takes value $1$ on each positive meridian
of $L$.
%(this class will be called {\it the orientation class} of the link and denoted by $\O$),
The class  $\O$ is fibered \ifff~ $L$ is an oriented fibered link.
For an oriented link $L$ the number  $\MNLo$ equals the usual Morse-Novikov number $\MM\NN(L)$.
The question of determining the invariant
$\MM\NN(L)$
was addressed in several papers.
In a series of papers
H. Goda gave a computation of Morse-Novikov numbers
$\MM\NN(K)$
for  prime knots $K$ with at most 10 crossings
\cite{Goda}, \cite{Goda2006bis} and
\mnn s
$\MM\NN(L)$
of prime
links $L$ with at most  9 crossings \cite{Goda2006}.
He proved in particular that for any prime link with at most 9 crossings its \mnn ~ is 2 or 0.
The invariant $\MM\NN(L)$ was first considered
in the paper \cite{PRW} of
A. Pajitnov, L. Rudolph and C. Weber,
where the basic properties of the invariant were studied and
the relations between  Novikov homology and
Alexander polynomial of links
were established.
%In the paper \cite{P-t} the third author proved that
%\bq\label{f:citing-tunn}
%\MM\NN(E(L),\Omega) \leq t(L)
%\end{equation}
%where $t(L)$ is the tunnel number of the link $L$.
Twisted Alexander polynomials and twisted Novikov homology
associated to representations of the fundamental group
allow to obtain an approximation
to the invariant $\MM\NN(L)$, see the papers of H. Goda, T. Kitano and T. Morifugi \cite{GodaKitanoMorifuji},
H. Goda and A. Pajitnov
\cite{GodaPajitnov}.
These invariants were discussed in the paper
 A. Pajitnov
\cite{Pajitnov-algan} for the case of arbitrary classes $\xi$.
%is less than or equal to twice the tunnel number of $L$.
The additivity of the number $\MM\NN(K)$
\wrt~ connected sum of knots was established in the article
F. Manjarrez-Gutierrez \cite{ManjarrezGutierrez} for the case of small knots,
and in the article of
K. Baker \cite{Baker} in general case.

%In the paper \cite{BakerManjarrezGutierrez} K. Baker and F. Manjarrez-Gutierrez
%introduced  an invariant $h(M,\g,\xi)$ for a sutured 3-manifold $(M,\g)$.

%-----------------
%In the present paper we develop
%algebraic and geometric tools, that allow us
%to compute the number
%$\MM\NN(E(L),\xi)$
%for arbitrary classes $\xi$ in many cases.

Returning to the case of general manifolds,
%As for upper bounds for $\MNxi$
we prove in the present  paper the following upper estimate:
\bq\label{f:tunn-m}
\MNxi\leq \MM(M;\pr M )-2
\end{equation}
\noindent
where $\MM(M;\pr M)$ is the usual Morse number of the manifold,
that is, the minimal number of critical points of a real-valued
Morse function $f$ on $M$ \sut~ $f~|~\pr M$ is constant.
In the case of complements to 2-knots in 4-sphere
an inequality equivalent to \rrf{f:tunn-m}
was used in the papers of H. Endo and A. Pajitnov
\cite{EndoPajitnov2017},
\cite{EndoPajitnov2017bis}
to determine Morse-Novikov numbers of 2-knots and surface-links.
For the case of exteriors of links this inequality
is equivalent to

\bq\label{f:tunn-L}
\MNLxi\leq 2t(L)
\end{equation}
where $t(L)$ is the tunnel number of $L$.

If $\xi$ is the orientation class for an oriented link $L$ in 3-sphere, the inequality
\rrf{f:tunn-L}
was proved in the paper \cite{P-t}. The proof of the general inequality
\rrf{f:tunn-m} repeats the argument of \cite{P-t}.
Combining the estimate \rrf{f:tunn-L}
with the fibering obstructions from Section 1 we determine the numbers
$\MNLxi$ for all 2-component prime links $L$ with $\leq 8$ crossings.
The results are summarized on  Figures 1 and 2
in Section \ref{s:6-7-8}.

In a recent paper
\cite{BakerManjarrezGutierrez}
K. Baker and F. Manjarrez-Guti\'errez
introduced and studied an invariant $h(M,\g,\xi)$
where $(M,\g)$ is a sutured 3-manifold and $\xi$ a 2-dimensional homology class. Presumably this invariant is related to Morse-Novikov numbers.

In the Section \ref{s:a_k}
we consider one special family of links.
Orient the annulus $S^1\times [-1,1]$ arbitrarily
and embed it in $S^3$ in such a way
that the core of the annulus be unknotted
and the linking coefficient of the boundary components
%(oriented as the boundary of the annulus)
be equal to $k$.
This link will be denoted by $a_k$.
For an $n$-uplet of natural numbers
$I=(k_1, ..., k_n)$ where $k_i\geq 1$
put
$$
a_I
=a_{k_1}\krest a_{k_2}\krest  ... \krest a_{k_n}.$$
We compute the \mnn s of the links $a_I$.
It turns out in particular that for some choice of orientation of the link $a_I$ the link is fibered, while
for another choice of orientation its \mnn~ is equal to $2n$.

In the Section \ref{s:irrat} we consider the generalizations
of some of the results of the previous sections to the case
of closed 1-forms with non-degenerate zeros.
The Section \ref{s:conj} contains some
 open questions.

\subsection*{Acknowledgements}
The first author was supported by JST SPRING, Japan Grant Number JPMJSP2180.
The second author thanks the Nantes University, the DefiMaths program, and ``Mission Invite'' for the support and warm hospitality.
The second author was supported by JSPS KAKENHI Grant Numbers 24K06707, 20K03578.
The third author was supported by the  JSPS International Fellowships for Research in Japan (Short-term), FY2022 Research Abroad and Invitation Program for International Collaboration, and the WRH program of Tokyo Institute of Science. He thanks Tokyo Institute of Science for the support and warm hospitality.
 \newpage

%\bibliography{mybib}
%\bibliographystyle{plain}
%\end{document}

\newpage
 \section{Novikov rings and Novikov homology}
 \label{s:nov-rings}

In this section we recall the basic facts from Novikov theory.
Let $G$ be a group; we denote by $\L$ the group ring $\zz G$.
Let $\xi:G\to\zz$ be a \ho.
Intuitively, the Novikov completion $\Lxi$
of the ring $\L$ consists of some special infinite linear combinations
of the elements of $G$, namely the linear combinations that are infinite
in the direction of decreasing of $\xi$. To give a precise definition
let $\lLL$ be the set of all formal linear combinations
(infinite in general)
$\l=\sum_{g\in G} n_g g, \ n_g\in \zz$.
For $\l\in \lLL$ put $\supp \l = \{g~|~ n_g\not=0\}$.
For $C\in\rr$ put
$$
[\l]_{\xi, C}
=
\{g\in \supp \l~|~  \xi(g)\geq C\}.$$

\bede\lb{d:l-xi}(\cite{NovikovDAN}, \cite{SikoravThese})
The set
\bq\lb{f:lxi}
\{\l\in\lLL~|~ \forall C {\rm \ \ the \ \ set \ \ } [\l]_{\xi, C} {\rm \ \ is \ \ finite \ \ }\}
\end{equation}
is called {\it Novikov ring}
and denoted by $\Lxi$.
%\bq\lb{f:lhat-min}
%\lLL_\xi^-=\{\l\in\lLL~|~ \supp \l \sbs \xi^{-1}(]-%\infty, 0]) \}
%\end{equation}
%\bq\lb{f:lxim}
%\Lximm = \Lxi \cap \lLL_\xi^-.
%\Lximm=\{\l\in\lLL~|~ \forall C {\rm \ \ the \ \ set \ \ } [\l]_{\xi, C} {\rm \ \ is \ \ finite \ \ }\}
%\end{equation}
\end{defi}

When $G=\zz$ and $\xi:\zz\to\zz$ is the identity \ho~
the ring $\Lxi$ is isomorphic to the Laurent power series ring in one variable with integer coefficients having only finite number of negative powers of $t$.

Let $X$ be a connected space, and $\xi\in H^1(X)$.
Let $\CCCC$ be a regular covering $\cov \ox p X$ of  $X$
with structure group $G$ \sut~ $p^*(\xi)=0$.
%and $\xi:G\to\zz$ a \ho.
The couple $(\CCCC,\xi)$ will be called {\it Novikov data} for $X$.
The class $\xi$ determines a \ho~ $G\to \zz$; by a certain abuse of notation we will denote this class by the same letter $\xi$.
%Let $\L=\zz G$ and $\Lxi$ be the Novikov completion of $\L$ \wrt~ $\xi$.
To a Novikov data $\EEEE = (\CCCC,\xi)$ we  associate the Novikov homology module
$$
\HN X \Lxi
=
H_*\big(\Lxi\tens{\L} S_*(\ox)\big),
$$
where $S_*(\ox)$ denotes the singular chain complex of $\ox$.
The Novikov homology corresponding to the universal covering $\wi X\to X$
%and the \ho~ $\pi_1(X)\to\zz$ induced by $\xi$
will be called {\it universal Novikov homology}.
The Novikov homology is related to critical
points of circle-valued Morse maps in the following way.

\beth[\cite{patoul}]
Let $M$ be a compact connected \ma, $\xi\in H^1(M) $ a regular
cohomology class, and $f:M\to S^1$ a $\smo$ map
belonging to $\xi$ such that $f~|~\pr M$ is
a fibration.
Let $(\CCCC,\xi )$ be  a Novikov data for $M$.
%such that $p^*([f])=0$.
Then there is a chain complex $\NN_*(f)$ of
free $\Lxi$-modules, freely
generated in degree $i$ by critical points of $f$ of index $i$, and such that
\bq\label{f:novik-iso}
H_*(\NN_*(f))
\approx
\HN X \Lxi.\end{equation}
\end{theo}

\bere
In \cite{patoul} the theorem was proved for closed \ma s.
The argument extends to the case of
non-empty fibered boundary without changes.
Observe that the covering $\ove{X}$ is not necessarily connected.
\enre

In particular,  the condition of vanishing of the universal Novikov homology
is necessary for the class $\xi$ to be fibered.
In case of 3-manifolds the inverse theorem is true,
it follows by combining a theorem of
J.-C. Sikorav and a theorem of  R. Bieri, W. Neumann and R. Strebel.
Their results imply the following:

\beth[\cite{BNS}, \cite{SikoravThese}]\label{t:bns-siko}
Let $M$ be a compact connected 3-manifold,
and $\xi\in H^1(M)$ be a regular class.
%Let $G=\pi_1(M),$  put $\L=\zz G$.
Then $\xi$ is fibered if and only if
the universal Novikov homology
$H_1(M, \Lxi)$
equals zero.
\end{theo}

%\newpage

\section{The case of 2-bridge links}
\label{s:two-br}

In this section
we use the Morse-Novikov theory to establish a fibration criterion
for a 2-bridge 2-component link  formulated in terms of
its Alexander polynomials.
%Let us  begin by recollection of some basic facts about Novikov
%rings.

\subsection{Fox calculus and universal Novikov homology}
\label{su:fox_nov}

We need some more definitions.

\bede\label{d:heights}

Let $\xi:G\to \zz$ be a \ho, $\L=\zz G$ and $ H=\Ker\xi$.
For an element
$x\in \Lxi, \ x= \sum_{g\in G} n_g g$
we denote by $h(x)$ the
number $\max \xi(g)$ where $g$ runs over
$\supp x$.
Set $
x_0
=
\sum_{\xi(g)=h(x)} n_g g.
$
%and write $x=x_0+x_1, \ h(x_1)<h(x_0)$.
We say that $x_0$ is the {\it principal part }
of $x$. An element $x$ is called $\xi$-monic, if
$x_0=\pm g, \ g\in G$.
\end{defi}

We have a simple Lemma (the proof will be omitted):
\bele\label{l:inv}
\been\item
Assume that $h(x)=0$, so that $x_0\in\zz H$. If $x_0$ is invertible in $\zz H$ then $x$ is invertible in $\Lxi$.
\item
If $x\in\Lxi$ is $\xi$-monic, then it is invertible.
\item
If $G\approx \zz^n$, then
$x\in\Lxi$ is invertible
\ifff~
it is $\xi$-monic. $\qs$
\enen
\enle
For a finitely presented group we can use Fox calculus to compute the Novikov homology. Recall that for a group $G$ with $m$ generators
$g_j$
and $n$ relations $r_i$
we have a chain complex $F_*$
of free left $\L$-modules
\bq
0\arl
C_0
\arrl \e
C_1
\arrl D
C_2
\arl
0
\end{equation}
with
$\rk\ C_0=1, \
\rk\ C_1=m,\
\rk\ C_2=n.
$
The matrices of the boundary operators are given by

$$
M(\e)=\begin{pmatrix}
       1-g_1, \ldots  , 1-g_m \\
      \end{pmatrix}, \ \
      M_{ij}(D) = \left( \dfox {r_i}{g_j}\right).
$$
where %$g_i$ are the generators, $r_j$ the relations, and
$\dfox {r_i}{g_j}$ are the Fox derivatives.
This chain complex will be called {\it the Fox complex}
associated to the presentation of $G$.

If $X$ is a connected finite CW complex with $\pi_1(X)\approx G$,
we have
a chain map $I:F_*\to C_*(\wi X)$ inducing an isomorphism in $H_0$ and $H_1$.
The next  lemma from homological algebra
enables us to use the Fox complex for computations with Novikov homology.
\bele\label{l:fox_novi}
Let $R$ be a $\L$-module.
The chain map  $I$ induces isomorphisms
$$
H_i(R\tens{\L} F_*)
\approx
H_i(R\tens{\L} C_*(\wi X))
$$
for $i=0,1$. $\qs $
\enle
In our applications the most important case is
$R=\Lxi$ (where $\xi:G\to\zz$
is a \ho). In this section we are interested in the case when $m=2, n=1$. For this case we have a simple sufficient condition for
vanishing of $H_1(R\tens{\L} F_*)$.

\bele\label{l:alg-lem-2}
Let
$$C_*=\{0\arl   C_0 \arrl {\pr_1} C_1\arrl  {\pr_2} C_2\}$$
be a chain complex of free left $R$-modules  over a ring $R$
with
$\rk\ C_0=1, \
\rk\ C_1=2,\
\rk\ C_2=1.
$
Write the matrices of the boundary operators as follows
$$
M(\pr_1)=
\begin{pmatrix}
       a_1,
       a_2
      \end{pmatrix},
      \ \ \
M(\pr_2)=
\begin{pmatrix}
       b_1 \\
       b_2
      \end{pmatrix},
      $$
Assume that  $a_2$ is invertible.
Consider two conditions:

(C1) \ \ $b_{1}$ is invertible.\ \ \ \ \ \ \ \
(C2) \ \
 $H_1(C_*)=0$.

 \noindent
Then
$(C1)\Rightarrow (C2)$ and if $R$ is commutative
then also $(C2)\Rightarrow (C1)$.
%\item
%If   $b_{1}$ is invertible, then
%% $H_1(C_*)=0$.
% \item If $R$ is commutative then the condition
%$H_1(C_*)=0$ implies invertibility of $b_{1}$.
% \enen
\enle

\Prf
To prove the first implication,
let $f_1, f_2$ be the free generators of $C_1$ and $g$
be the free generator of $C_2$.
Then  $\Ker \pr_1$ is a free submodule of $C_1$ generated by
$p=f_1-a_1a_2^{-1} f_2$. To check that $p\in \Im \pr_2$, put
$q=b_{1}^{-1}g_j$.Then
$\pr_2(q)=
b_{1}^{-1}b_{2}f_2+f_1$.
The condition $\pr_1\circ\pr_2=0$
implies
$b_{1}^{-1}b_{2} = -a_1a_2^{-1}$,
so $\pr_2(q)=p$ and the proof of the first implication is over.
The proof of the second one  is left to the reader. $\qs$

%\[
%%        \begin{bmatrix}
 %       0 & 1 & 0 & 0\\
 %       0 & 0 & 1 & 0\\
 %       0 & 0 & 0 & 1\\
 %       -1 & 9 & -14 & 8\\
 %       \end{bmatrix}.
 %   \]

 Together with Theorem \ref{t:bns-siko}
this leads to the following sufficient condition of fibering for
\tbtc s and generalized orientations.

\bepr\label{p:tbtc-fibr}
Let $L$ be a \tbtc. Choose a presentation of $\pi_1(E(L))$
with two generators $a,b$ and one relation $r$.
Let $\xi\in H^1(E(L))$ \sut~ $\xi(b)\not=0$.
Assume that $\prpr  r w $
is $\xi$-monic. Then $\xi$ is fibered.
\enpr
\Prf
The element $1-b$ is invertible in $\Lxi$,
so by
Lemma \ref{l:alg-lem-2}
we have
$H_1(E(L), \Lxi)=0$ and Theorem \ref{t:bns-siko}
implies the required assertion. $\qs$

%\newpage

\subsection{Abelian coverings and Alexander polynomials}
\label{su:abel_alex}

Let $L$ be an $n$-component link in $S^3$,
denote $\pi_1(E(L))$ by $G$ and $H_1(E(L))$ by $H$,
so that $H$ is a free abelian group
of rank $n$.
The group rings of these groups will be denoted by $\L$, respectively $\PPPP$,
so that $\PPPP$ is a Laurent polynomial ring.
The Hurewicz \ho~ $p:G\to H$ extends to a ring morphism
$\wi p: \L\to\PPPP$.
Let $\UUUU:\wi {E(L)}\to E(L)$ be the universal covering,
and
$\AAAA:A(L)\to E(L)$
be the maximal abelian
covering. The structure group of $\AAAA$
is isomorphic to $H$.
%is a free abelian group $H$ of rank $n$,
%isomorphic to $H_1(E(L))$.
Let $\xi:G\to\zz$ be a \ho.
It factors through a \ho~
$\eta:H\to\zz$.
We have corresponding Novikov completions
$\Lxi$ and $\Peta$ and a ring morphism
$\Lxi\to\Peta$, extending $\wi p$.
We have an isomorphism
$$
\Peta\tens{\PPPP} C_*(A(L))
\approx
\Peta\tens{\Lxi}
\left(\Lxi\tens{\L} C_*(\wi {E(L)})\right).
$$

Therefore if the universal Novikov homology equals zero,
the same is true for the abelian Novikov homology
$H_*(E(L),\Peta)$. This provides a necessary condition
for fibering, which in turn can be reformulated in terms of Alexander
polynomials. Let us recall the necessary definitions.
Write
$$
C_*(A(L))=\{ 0\arl  C_0\arrl {\pr_0} C_1\arrl {\pr_1} C_2 \arrl {\pr_2}\ldots \}
$$
This is a free chain complex over the Laurent polynomial ring $\PPPP\approx \zz[\zz^n]$.
Observe that  $\rk\ C_0=1$, and denote $\rk\ C_1$ by $s$.
\bede
The gcd of $(s-1)$-minors of the matrix of $\partial_2$
is called the {\it Alexander polynomial}
of $L$, and denoted by $\D(L)$.
\end{defi}

One can show that this polynomial is a chain homotopy invariant of the chain complex
$C_*(A(L))$, see for example\cite{PajAlgan}.
\footnote{
The polynomial $\D(L)$ is defined up to multiplication by an
invertible element of the ring $\PPPP$, and in the sequel
the equalities containing the polynomial $\D(L)$
must be understood modulo thie indeterminacy.}

The next proposition is proved in
a more general setting in \cite{PajAlgan}, Theorem 5.5.
\bepr\label{p:monicAlex}
Let $\xi:\zz^n\to\zz$ be a \ho.
The two following properties are equivalent:
\been\item
$H_1\left(\Peta\tens{\PPPP}C_*\left(A(L)\right)\right)=0$,
\item
$\D(L)$ is $\eta$-monic.
\enen
Therefore if the class $\xi$ is fibered, then
the polynomial $\D(L)$ is $\eta$-monic.
\enpr

In the case of \tbtc~ the Alexander polynomial can be explicitly computed
in terms of the Fox derivatives.
%Let $G$ denote the fundamental group of $E(L)$,
%and let $p:G\to H$ be the Hurewicz \ho.
%Denote by $\wi p:\zz[G]\to\zz[H]$ the corresponding \ho~
%of group rings.
Choose a presentation of $\pi_1(E(L))$
with two generators $a,b$ and one relation $r$.
Let $t_1=\wi p(a),
t_2=\wi p(b).$
The proof of the next Proposition is obvious.
%Denote by $t_1$ the image of $a$, and by $t_2$ the image of $b$ in $H$.
%The elements $t_1, t_2$ can be considered
%as the elements of the group ring $\zz[H}$.
%Denote by $\wi p:\zz[G]\to\zz[H]$ the \ho~
\bepr\label{p:alex_2br}
The element
$D=\wi p (\prpr r b ) $
is divisible by $1-t_1$
and we have
\bq\label{f:alex_2br}
\D(L) = D/(1-t_1).
\end{equation}
\enpr

%\newpage

The main aim of this section is the next Theorem.
\beth\label{t:2br}
Let $L$ be a \tbtc, and $\xi\in H^1(E(L))$ a regular cohomology
class.
%\been\item[(A)]
The following two conditions are equivalent:

(1) The class $\xi$ is fibered. \ \ \
(2) $\D(L)$ is $\xi$-monic.
%item[(B)]
%If $m,n\in\zz$ are odd and $\gcd(m,n)=1$ then the two conditions above
%%%3) $d(L,\xi)$ is monic.
%\enen
\enth
%Observe that the point (B) implies the following result:
%\begin{quote}
% {\it A \tbtc~ is fibered \ifff~
% its Alexander polynomial is monic.}
% \end{quote}
%This is a particular case of a K. Murasugi's theorem
%\cite{Murasugi}, saying that an alternate link is fibered
%\ifff~ its Alexander polynomial is monic.

%Observe that theorem \ref{t:2br} (B) generalizes  the classical result
%of K. Murasugi \cite{Murasugi}, saying that
%a \tbtc~ link is fibered if its Alexander polynomial is monic.

The implication (1)$\Rightarrow$(2)
is already proved.
The proof of the inverse implication occupies the rest of the section.

\subsection{$p$-irreducible elements in group rings}
\label{su:p-irred}

This subsection is purely algebraic.
We prove here some auxiliary results to be used in the next subsection.
Let $G,H$ be any groups an $p:G\to H$ a group \ho.
Extend it to a ring \ho~ $\wi p : \zz G \to \zz H$.
Let $\l\in\zz G, \ \l=\sum_{g\in G} n_g g$.
\bede\label{d:irred}
We say that $\l$ is {\it $p$-irreducible}
if for any $g\in \supp\l$ the sign of $n_g$ is determined
by $p(g)$ (that is, all numbers $n_h$ with $h\in p^{-1}(p(g))$ have the same sign).
\end{defi}

The $p$-irreducible elements $\l$ have the following
properties:
\been\item
$p(\supp \l)=\supp (\wi p (\l)$).
\item Write $\wi p (\l)=\sum_{h\in H} m_u u$. Then
$|m_u|= \sum_{g\in p^{-1}(u)} |n_g|$.
\enen

Let $\eta:H\to\zz$ and $\xi=\eta\circ p$.
Let $\l$ be a $p$-irreducible element of $G$, and denote by $\l_0$
the $\xi$-principal part of $\l$. The properties 1) and 2) above
imply the following:

\been\item[3)]
The $\eta$-principal part of $p(\l)$ equals $\wi p(\l_0)$.
\item[4)] \label{f:monic}
$\wi p(\l)$ is $\eta$-monic \ifff~ $\l$ is $\xi$-monic.
\enen

\bere\label{r:two-groups}
Observe that if we have morphisms
$G_1\arrr p G_2 \arrr q H$, where $q$ is surjective, and $\l\in\zz G_1$ is $pq$-irreductible, then
$p(\l)\in \zz G_2$ is $q$-irreductible.
\end{rema}

\begin{examples}\label{e:examples}  1)
Let $G$ be  a free group on 2 generators $a,b$. Put $H=G/[G,G]$
and let $p:G\to H$ be the projection.
Denote by $d:G\to\zz$ the \ho~ defined by $d(a)=d(b)=1$.
We say that $\l\in\zz G ,
\l=
\sum_{g\in G} n_g g
$ has
{\it alternating coefficients}
if
$
\sgn (n_g) = (-1)^{d(g)}.
$
An element having alternating coefficients
will be also called {\it $AC$-element}.
It is clear that any $AC$-element is $p$-irreducible, since
the \ho~ $d$ factors through $H$.

2)  Let $g$ be a free group on 2 generators $a,b$, and $H=\zz$, and
$p:G\to H$ is defined by $p(a)=m,\ p(b)=n$.
Assume that both $m,n$ are odd. Then any $AC$-element in $G$ is
$p$-irreducible.
\end{examples}

A particular class of examples of $AC$-elements is provided
by Fox derivatives of elements of special type.

\bede\label{d:hartl}
A word $w$ in the alphabet of two letters $a,b$ is called
{\it a word of type (H)}
if it equals $1$ or is of the form
$$
w= b^{\e_1}a^{\eta_1} \cdots $$
where $\e_i, \eta_i = \pm 1$.
(That is, the word begins by $b$ or $b^{-1}$
and then the letters $a$ and $b$ appear in turn.)
\end{defi}

\bepr\label{p:typeH}
If $w$ is of type (H), then every monomial in $\dfox w b $
is of type (H), and the element  $\dfox w b \in \zz G$
has alternating coefficients. \enpr
\Prf
 Induction in the length $l(w)$ of the word $w$.
For $l(w)=1$ we have $w=1$ or $w=b$ or $w=b^{-1}$.
Then
$\dfox w b = 0$ or
$\dfox w b = 1$ or $\dfox w b = -b^{-1}$.

Assuming that the assertion holds for every word of length $<k$, consider
a word $w$ with $l(w)=k$. Then either $w=b a^\nu u$ or
$w=b^{-1} a^\nu u$ with $\nu=\pm 1$ and $u$ is a word of type (H)
of length $\leq k-2$.
We have
$$
\dfox w b = 1 + b a^\nu \dfox u b
\ \ \ \
{\rm or ~ }
\ \ \
\dfox w b = -b^{-1} +  b^{-1}a^\nu \dfox u b.
$$

%In the second case we have
%$\dfox w b = -b^{-1} +  b^{-1}a^\nu \dfox u b$.
\noindent
Apply the induction assumption and the proof is over.
 $\qs$

\subsection{Proof of Theorem \ref{t:2br}}
\label{su:proof2br}

The proof relies on the special presentation of the fundamental group
$G$
of \tbtc~ links, see \cite{Hartley}, \cite{Hoste}.
Namely, the group has two generators, $a,b$ and one relation

\bq\lb{f:pres}
r=awa^{-1}w^{-1}, \
{\rm and ~ }
w=
b^{\e_1}a^{\e_2}\cdots b^{\e_{q-1}} \
{\rm with ~ }
\e_i=(-1){\lfloor ip/q\rfloor}.
\end{equation}
Let $F_2$ be the free group on two generators $a,b$;
observe that the word $w$ is of type $(H)$.
We have
$$
\dfox r a
=
1-w+(a-1)\dfox w a,
\ \ \ \ \
\dfox r b
=
(a-1)\dfox w b.
$$

The class $\xi$ being regular, one of the elements
$\xi(a), \xi(b)$ is non-zero. assume for example that
$\xi(a)\not=0$.
According to Proposition \ref{p:typeH}
the element $\dfox w b$ is  of $AC$-type,
and is therfore $p$-reducible, where $p:F_2\to H$
is the projection of  $F_2$
 onto its first homology group $H$, isomorphic to $\zz^2$.
The map $p:F_2\to H$ factors as
$F_2\to G \arrr q  \zz^2$ and applying Remark \ref{r:two-groups}
we see that that element
$\dfox w b $ in the group $G$ is $q$-irreducible.
The element $\wi p(\dfox w b )$ is $\eta$-monic,
therefore  $\dfox w b $
is $\xi$-monic
(see the property 4), page  \pageref{f:monic}).
Therefore
according to Proposition \ref{p:tbtc-fibr} the
class $\xi$ is fibered. $\qs$

\subsection{Fibered classes and 1-variable Alexander polynomials}
\label{su:1var}

In some cases the abelian invariants deduced from infinite cyclic
coverings are sufficient to detect fibred classes in cohomology
of \tbtc~ links.

Let $L$ be a \tbtc~ link, put $G=\pi_1(E(L))$, let $\xi:G\to \zz$ be an epimorphism.
We have the corresponding connected abelian covering $B(L))\to E(L))$
with structure group $\zz$. Denote its group ring by $P$.
The chain complex of $B(L)$ is of the form
$$
C_*(B(L))=\{ 0\arl  C_0\arrl {\pr_0} C_1\arrl {\pr_1} C_2 \arrl {\pr_2}\ldots \}
$$
where $C_i$ are free modules over $P$, and $\rk C_0=1$.
Denote $\rk\ C_1$ by $r$.

\bede
The gcd of $(r-1)$-minors of the matrix of $\partial_2$
is called the {\it Alexander polynomial of $L$ \wrt~ $\xi$}
 and denoted by $d(L,\xi)$.
\end{defi}
We will use the presentation \rrf{f:pres}
of the group $G$. The \ho~ $\xi$ is determined by its values on $a$ and $b$.
Put $\xi(a)=m, \xi(b)=n$, surjectivity of $\xi$ mplies that $\gcd(m,n)=1$

\beth\label{t:alex1var}
Assume that $\gcd(m,n)=1$ and both $m,n$ are odd.
Then the class $\xi$ is fibred \ifff~
$d(L,\xi)$ is monic (that is, its  coefficient of its principal term equals $\pm 1$)
\enth

The proof repeats the proof of Theorem \ref{t:2br}
(using Example 2, page \pageref{e:examples} instead of Example 1)
We omit the details. $\qs$

Observe that this theorem implies a particular case of K. Murasugi's theorem about fibering of
alternate links, namely we have the following Corollary

\beco [\cite{Murasugi}]
Let $L$ be a \tbtc. Then $L$ is fibered \ifff~~
its Alexander polynomial is monic.
$\qs$
\enco

\section{Connected sums of links, and Gabai's theorem}
\label{s:conn-sum}

We begin by  preliminaries about Novikov homology
of induced coverings.
Let $X\sbs Y$ be finite CW-complexes. Put
$$G=\pi_1(Y), \ \  \L=\zz G, \ \
H=\pi_1(X), \ \ \G=\zz H
$$

Let $\xi:G\to\zz$ be a \ho~ and denote by $\eta:H\to\zz$ the induced \ho.
We denote by $\wi X, \wi Y$ the universal coverings and by $\ove X\to X$
the covering induced from the universal covering $\wi Y \to Y$.
For $X$ we have two Novikov datas:
the one corresponding to the universal covering with Novikov
homology $\HN { X } \Geta$ and the one correspondng to the
covering $\ove X\to X$ with Novikov
homology $\HN { X } \Lxi$
We have natural \ho s
$$
\HN {X} \Lxi \arrr {i}\HN { Y } \Lxi,
\ \ \
\HN { X } \Geta\arrr {j}\HN { Y } \Lxi.
$$

\bepr\label{p:monos}
Assume that
$i$ is injective, and moreover, that
 $\pi_1(X) \to \pi_1(Y)$ is injective.
 Then  $j$ is injective.
 \enpr
 \Prf
  Choose a connected component
  ${\ove X}_0$
  of $\ove X$.
  The map
  ${\ove X}_0
  \to
  X$
  is then a regular covering with structure group $K=\{g\in G ~|~ g({\ove X}_0)={\ove X}_0\}$.
Observe that
    ${\ove X}_0$
    is 1-connected. Indeed
    the \ho
    $$
    \pi_1({\ove X}_0) \to \pi_1(X)\to \pi_1(Y)$$
    is injective by assumption, and it factors through the group $\pi_1(\wi Y)$ which
    is trivial.
    Therefore the map ${\ove X}_0\to X$ is the universal covering, we have $H\approx \pi_1({\ove X}_0)$,
    and it suffices to prove the next lemma.

\bele\label{t:j-mono}
The homomorphism
$$
\HN {{ X}} \Geta \arrr {J}\HN {X} \Lxi
$$
induced by inclusion ${\ove X}_0 \sbs \ox$, is injective.
\enle

\Prf
In case when $X$ is a finite CW complex there is an alternative description of the Novikov homology using cellular chains. Choose a continuous function $f:X\to S^1$ belonging to $\xi$
and choose a $\xi$-equivariant lift $F:\ox\to\rr$ of $f$
(that is, $F(gx)=F(x)+\xi(g)$).
For each cell $\s$ of the cell decomposition of $\ox$ denote by $h(\s)$ the minimum  of $F$ on $\s$.
A formal linear combination
$\l=\sum_{i=1}^\infty n_i\s_i$
where $n_i\in\zz$ and $\s_i$ are $k$-cells of $\ox$,
will be  called {\it regular}
if for any $C\in\rr$ there is only a finite set of $\s\in supp(\l)$
with $h(\s)>C$.
The set of all regular linear combinations of $k$-cells will be denoted by
$T_k(\ox)$, this is a free $\Lxi$-module.
We have then
$$
H_*\big(T_*(\ox)\big)\approx \HN  X \Lxi.
$$

To construct a left inverse $Q$  for $J$,
define $Q$ on the cells of $\ox$ as follows:
$Q(\s)=\s$ if $\s$ is a cell of $\oox$, and $Q(\s)=0$ if not.
Extend $Q$  to $T_*(\ox)$ in the obvious way.
Then $Q$ commutes with boundaries
(since for a cell of $\oox$ its boundary is also in $\oox$)
and we have obviously $Q\circ J=\Id$.  $\qs$

\bere
The map $Q$ is a \ho~of abelian groups but in general it does not commute with the action of Novikov rings.
\enre

Let $L_1, L_2$ be oriented links in $S^3$.  Choose a component
$K_1$ of $L_1$ and a component
$K_2$ of $L_2$. In \cite{Hashizume}
the connected sum $L=L_1\krest L_2$ of links $L_1$ and $L_2$
along the components $K_i$ was defined.
We will use a $\smo$ version of this definition. The link
$L_1\krest L_2$
has the following properties.

\been\item
The sphere $S^3$ is the union of two 3-dimensional disks $D_1$ and
$D_2$ intersecting by a sphere $\Sigma$, the link $L$ is
transverse to $\Sigma$ and $L\cap\Sigma=\{\a, \b\}, \ \
L=(L_1\cap D_1)\cup (L_2\cap D_2).
$
\item
There are two arcs $A_1$ and $A_2$ in $K_1$, resp. $K_2$, \sut~
$$
L_1\cap D_1=L_1\sm \Int A_1, \ \  \ L_1\cap D_2 = A_1
$$
$$
L_2\cap D_2=L_2\sm \Int A_2, \ \ \
L_2\cap D_1
= A_2.
$$
\enen

Let $R= D_1\sm L_1,  \ S= D_2\sm L_2$.
Denote the twice punctured sphere
$\Sigma \sm \{\a, \b\}$ by $\Sigma_0$.
Then we have
$$S^3\sm L = R\cup S,\ \ R\cap S= \Sigma_0.$$
Observe that $R$ has the homotopy type of
$S^3\sm L_1$ and
$S$ has the homotopy type of
$S^3\sm L_2$.
Consider the universal  Novikov homology
$H_1(S^3\sm L, \Lxi)$
where $\L=\zz [\pi_1(S^3\sm L)]$
and $\xi\in H^1(S^3\sm L)$.
Denote by $\xi_i$ the image of $\xi$ in $H^1(D_i\sm L_i)$ (here $i=1,2$. )

\beth\label{t:links}
Let $\xi\in H^1(C_L)$ be a regular class \sut~
the restriction of $\xi$ to $\Sigma_0$ is non-zero.
Then the two following conditions are equivalent
\end{theo}
\been
\item$\xi$ is fibered.
\item both $\xi_1, \xi_2$ are fibered.
\enen
\Prf
Only the implication $(2)\Rightarrow (1)$ has to be proved.
Denote by
$p:\wi{S^3\sm L} \to S^3\sm L$ the universal
covering. We have
$$
\wi{S^3\sm L}
=
p^{-1}(R)
\cup
p^{-1}(S),$$
and
$$
p^{-1}(R)
\cap
p^{-1}(S)
=
p^{-1}(\Sigma_0).
$$
Consider the corresponding Mayer-Vietoris sequence
for the homology
$H_1(S^3\sm L, \Lxi)$.
The restriction
of $\xi$ to $\Sigma_0$ is not zero, so we have
$$
H_*(
R\cap S, \Lxi)
\approx
H_*(\Sigma_0,\Lxi)=
0,
$$
and the Mayer-Vietoris exact sequence
implies an  isomorphism
$$
H_*(S^3\sm L, \Lxi)
\approx
H_*(R, \Lxi)
\oplus
H_*(S, \Lxi).
$$
Therefore the inclusions $R\sbs S^3\sm L,\ S\sbs S^3\sm L$
induce monomorphisms in the Novikov homology with
coefficients in $\Lxi$.
Both inclusions induce monomorphisms in $\pi_1$.
The previous Proposition
implies that the universal Novikov homology of
$R$ and $S$ vanishes, so it is also the case for
$S^3\sm L_1$ and $S^3\sm L_2$ and the proof is over. $\qs$

A particular case of this theorem implies the  following theorem of Gabai.
\beco\cite{Gabai}
Let $L_1, L_2$ be oriented links in $S^3$.
Then $L_1\krest L_2$ is fibered if and only if
both $L_1$ and $L_2$ are fibered. $\qs$
\end{coro}.

%\newpage

\section{Morse-Novikov numbers of manifolds}
\label{s:mn-manif}

Let $M$ be a compact manifold with $\pr M\not=\ems$.
The pair $(M; \pr M)$ will be considered as a cobordism between $\pr M$ and $\ems$.
\bede\label{d:m-func}
{\it A Morse function $f$ on $(M;\pr M)$ }
is  a $\smo$ function $f:M\to\rr$ \sut~
\been\item$f~|~\pr M$ is constant and $f(\pr M)$ equals the minimum of $f$.
\item
all critical points of $f$ are non-degenerate and belong to $M\sm \pr M$.
\enen
\end{defi}

For a Morse function $f$
we denote by
$m_i(f)$ the number of critical points of $f$ of index $i$,  by $S(f)$ the set of
 all critical points of $f$ and by $m(f)$ the cardinality of $S(f)$.

 \bede\label{d:m-numb}
 The minimal possible number of critical points of a Morse function on $(M; \pr M)$ will be called {\it the Morse  number of $(M; \pr M)$ } and denoted by $\MM(M; \pr M)$.
 A Morse function $f$ is called {\it minimal }
 if $m(f)=\MM(M; \pr M)$
\end{defi}

\bere\label{r:zero}
A standard cancellation argument implies that
if $M$ is connected then for a minimal Morse function $f$ we have
$m_0(f)=0, \ m_n(f)=1$, where $n=\dim M$.
\enre

We proceed now to Morse-Novikov numbers.
For a map $f:X\to S^1$  we denote by $[f]\in H^1(X)\approx [X, S^1]$
the cohomology class $f^*(i)$ where $i\in  H^1(S^1)$
is the fundamental class of $S^1$.
We say that $f$ {\it belongs to a cohomology class $\xi\in H^1(X)$}
if $[f]=\xi$.
Recall that a $\smo$ map $f:M\to S^1$ from a closed manifold $M$ to $S^1$
is called {\it Morse map} if all critical points of $f$ are non-degenerate.
\bede\label{d:regular}
For a compact \ma~  $M$ with  a non-empty boundary we  say
 that a $\smo$ map $f:M\to S^1$
is a  {\it Morse map} if all the critical points of $f$ are non-degenerate
and belong to  $M\sm \pr M$,
and the map $f~|~ \pr M:\pr M\to S^1$ has no critical points
(i.e. $f~|~ \pr M$ is a $\smo$ fibration).
\end{defi}
Observe that this condition implies that
$
[f~|~\pr M]\not=0$ and
$
[f]\not=0$. Recall that  that a class $\xi\in H^1(M)$
is called {\it regular} if there is a fibration $\phi:\pr M\to S^1$
\sut~ $[\phi]$ equals the image of $\xi$ in $H^1(\pr M)$.
In case of $\pr M=\ems$ any class will be called regular.
A standard Morse-theoretical argument
implies that for any regular class $\xi$ there is a Morse map
$f:M\to S^1$ belonging to $\xi$.
% %(see \cite{milnorhcob} for the case of real-valued %Morse maps)

\bede
Let $M$ be a  connected compact \ma.
Let  $\xi\in H^1(M)$ be a regular class. The {\it Morse-Novikov number $\MM\NN(M,\xi)$}
is the minimal possible number of critical points
of a Morse map $f:M\to S^1$ belonging to $\xi$. We say that a class $\xi\in H^1(M)$
is {\it fibered} if $\MM\NN(M,\xi)=0$.
\end{defi}

\bere\label{r:handle-num}
In the paper \cite{GodaSakasai}
H. Goda and T. Sakasai introduce the notion of
a {\it handle number }
$h(M,\g)$ of a 3-dimensional sutured manifold
$(M,\g)$.
%This invariant is closely related to
%\mnn s.
In the paper
\cite{BakerManjarrezGutierrez}
K. Baker and F. Manjarrez-Guti\'errez
define an invariant
$h(M,\g, \xi)$ where
$(M,\g)$ is a 3-dimensional sutured manifold
and $\xi\in H_2(M, \pr M)$.
These invariants are related to the Morse-Novikov numbers. Presumably we have
$h(E(L),\g_0,  \xi)=\MM\NN(L, \wh\xi)$
where $\g_0$ is the suture containing
the whole of $\pr E(L)$ and $\wh\xi$ is the 1-cohomology class dual to $\xi$.
\end{rema}

\subsection{Morse numbers versus Morse-Novikov numbers}
\label{su:m-vers-mn}

A Morse function on a compact manifold has necessarily a maximum.
This is not the case for a Morse map to a circle belonging to
a non-zero cohomology class.

\beth\label{t:m-vers-mn}
Let
$M$ be a connected compact manifold with $\pr M\not=\ems$,
and $\xi$ be a regular class in $H^1(M,\zz)$.
Then
$$
\MM\NN(M,\xi)\leq\MM(M; \pr M)-2.
$$
\enth
\Prf
%Choose a collar
%$\psi: \pr M\times [0,\e[ \to M$ for $\pr M$ in $M$.
Pick a minimal Morse function $f$ on $(M;\pr M)$.
Pick any $ \smo$ function $g:M\to S^1$ belonging to $\xi$
such that $g~|~\pr M$ is a fibration.
We can assume that
the function $g$ is constant on some \nei~ $U$ of $S(f)$.
Consider the closed 1-form $\o_C=dg+Cdf$ where $C>0$.
Observe that $\o$ equals $Cdf$ in $U$, and outside $U$ the minimum of $\grad (Cdf)$ becomes
arbitrarily large with $C$. Therefore
for $C$ sufficiently large the
form $\o_C$ has no zeeros in $M\sm U$
and we have
$$
m_i(\o_C)=m_i(f)
$$
Fix such a number $C$.
The cohomology class of $\o_C$ equals
that of $dg$, so that $\o_C=dh$ with $h:E(L)\to S^1$
a Morse map belonging to $\xi$, and with $m_i(h)=m_i(f)$
for every $i$.
The  Morse map $h$ has one critical point of index $n=\dim M$, and this point can be cancelled. The proof of the next two lemmas is done by standard rearrangement and cancellation arguments, similarly to
Lemmas 3.1 and 3.2 of \cite{PRW}, and will be omitted.

\bele\label{l:connect}
There is a regular Morse map
$k:M\to S^1$ in the class $\xi$ \sut~
$m_i(k)\leq
m_i(h)$ for every $i$ and one of regular
level submanifolds $k^{-1}(\l)$ of $k$
is connected.
\enle

\bele\label{l:cancell}
Assume that
$k:M\to S^1$
is a regular Morse map
in the class $\xi$ such that
one of its regular level submanifolds is connected.
Then there is a regular Morse map
$r:M\to S^1$ in the class $\xi$ \sut~
\bq\label{f:indices}
m_i(r)=m_i(k) \ \  {\rm for } \ \ i\leq n-2,
\ \
m_{n-1}(r)\leq m_{n-1}(k),
\ \
{\rm and  }\ \
m_n(r)=0.
\end{equation}
\enle

Now we can prove our theorem.
Apply Lemma \ref{l:connect} and consider the map $k$.
We have two possibilities:
(1) $m_i(k)=m_i(h)$ for every $i$,
and (2) for some $i$ we have $m_i(k)< m_i(h)$.
In the case (2) we have $m(k)\leq m(h)-2$
(since $m(k)\equiv   m(h)\equiv \chi(E(L)) \mod 2 $)
and the proof is over.
In the case (1) we apply Lemma \ref{l:cancell}
and obtain q regular Morse map $r$ with
$m_i(r)=m_i(h)$ for $i\leq n-2$,
amd $m_{n-1}(r)\leq m_{n-1}(h)$, and $m_n(r)=0$.
Counting Euler characteristic again we deduce that
$m_{n-1}(r)=m_{n-1}(h)-1$, therefore
$m(r)=m(h)-2$, and the proof is over also in thic case. $\qs$

\subsection{Morse-Novikov numbers and tunnel numbers of links in $S^3$}
\label{su:mn-links}

Recall the definition of the tunnel number of a link
$L\sbs S^3$.
An arc $A$ is $S^3$ is called {\it a tunnel of $L$}
if it intersects $L$ by two points that are
extremities of $A$.
A finite collection $
\AA=\{A_1, \ldots , A_n\}$
of tunnels is called {\it tunnel system for $L$, }
if $A_i$ are mutually disjoint.
Denote by $|\AA|$ the union of all $A_i$, by $T(A)$
an open tubular \nei~ of $L\cup |\AA|$, and by
$E(\AA)$ the set $S^3\sm T(A)$.

\bede\label{d:tun-num}(B. Clark, \cite{Clark})
Tunnel number $t(L)$ of $L$ is the minimal cardinality of a
tunnel system $\AA$, such that $E(\AA)$ is a handlebody.
%Tunnel number of $L$ will be denoted by $t(L)$.
\end{defi}

%\newcommand{\ClarkHG}
%{B. Clark,
%\emph{
%The Heegaard genus of manifolds obtained
%by surgery on links and knots},
%Internat. J.

We will use an equivalent  definition
which is more convenient in the  Morse-theoretical setting.

\bede\label{d:tunn-morse}
Tunnel number of $L$ is the
minimal possible
number of critical points of index $1$ of a Morse
function $f$
on $(E(L); \pr E(L))$.
\end{defi}

Observe that if a Morse function $f$ on $(E(L); \pr E(L))$
has minimal possible number $m_1(f)$ then
$m_0(f)=0$, and we can also assume that $m_3(f)=1$.
Given that $\chi(E(L))=0$ we deduce that
\bq\label{f:morse-tunnel}
\MM(E(L); \pr E(L))=2t(L)+2.
\end{equation}

\beth\label{t:tunn-ineq}
Let $L$ be a link in $S^3$ and $\xi\in H^1(E(L))$
a regular cohomology class. Then
$$
\MM\NN(L,\xi)\leq 2t(L).
$$
\enth
\Prf
Let $f$ be a minimal Morse function on $(E(L); \pr E(L))$.
%We have then $m_0(f)0, m_3(f)=1$.
%Observe that $\chi(E(L))=0= -m_1(f)+m_2(f)-1$,
%therefore
%$$
%m(f)=m_1(f)+m_2(f)+1 = 2m_1(f)+2 \geq 2t(L)+2.
%$$
Applying Theorem \ref{t:m-vers-mn}
we deduce that
$$
\MM\NN(L,\xi)
\leq \MM(E(L); \pr E(L)) -2
=2t(L). \ \ \ \ \qs
$$

\section{Links $a_k$ and their connected sums}
\label{s:a_k}

In this section we consider an infinite family of
links, for which the Morse-Novikov number
varies from $0$ to arbitrary large values.

We orient the annulus $S^1\times [-1,1]$ arbitrarily
and embed it in $S^3$ in such a way
that the core of the annulus be unknotted
and the linking coefficient of the boundary components
%(oriented as the boundary of the annulus)
be equal to $k$.
This link will be denoted by $a_k$.
Denote by $a_k^-$ this link endowed with
the orientation induced from the boundary of the annulus.
Reversing the orientation on one of the two components f the link,
we obtain a  link  which we denote $a_k^+$, this is a torus link.
\beth
Denote by $\xi_0$
the meridional class of the link $a_k^-$.

1) We have $\MM\NN(a_k,\xi_0)=2$.

2) Every regular cohomology class $\xi\in H^1(S^3\sm a_k)$ which is not
an integer  multiple of $\xi_0$ is fibered.
\enth
\Prf
Denote by $c$ and $b$ the meridional classes of $a_k^-$, and
set $ v=b, \ u = cb^{-1}$.
Using the standard Wirtinger presentation of the group
$
\pi_1(S^3\sm a_k)
$
it is easy to show that
\bq\lb{f:press}
\pi_1(S^3\sm a_k)
\approx
\langle u,v ~|~ [u^k, v]=1\rangle,
\end{equation}
\bq\lb{f:foxx}
\ruru=(1-v)(\ukk), \ \ \ \
\rvrv = u^k-1.
\end{equation}
Assume that the class $\xi$ is not an integer multiple of $\xi_0$.
Then  $\xi(u)\not=0$.
Therefore the elements
$1-u^k$ and $\ukk$ of the Novikov ring
$\Lxi$ are invertible,
%in the Novikov ring,
and applying
Lemma \ref{l:alg-lem-2} and Theorem
\ref{t:bns-siko}
we deduce the point 2) of our assertion.

Proceeding to point 1), observe that the link $a_k$ is a 2-bridge link, therefore $t(a_k)=1$ and $\MM\NN(a_k,\xi_0)\leq 2$.
In order to obtain the equality, we will use the
Novikov homology obtained from the maximal abelian covering
$\ove{E(a_k))} \to E(a_k)$.
The group ring $\AAAA$ of the structure group of
this covering is isomorphic to $\zz[u^{\pm 1},v^{\pm 1}   ]$.

Using\rrf{f:foxx} it is easy to show that
$$
H_1(E(a_k),\Axi)
\approx
\Axi /(P_k),
\ \
{\rm where } \ \
P_k=(\ukk).$$
Returning to the meridional variables
$c,b$ we have also

\bq\lb{f:qk}
H_1(E(a_k),\Pxi)
\approx
\Bxi /q_k, \ \ {\rm with } \ \
\BBBB= \zz[b^{\pm 1},c^{\pm 1 }],
{\rm and } \ \
\end{equation}

\bq\lb{f:qkk}
q_k= b^{k-1}+ b^{k-2}c + \ldots + bc^{k-2} + c^{k-1}.
\end{equation}

This homology module is not equal to zero, therefore
the class $\xi_0$ is not fibered. $\qs$

Let us proceed to connected sum of links $a_k$.

\bede
Let $I=(k_1, ..., k_n)$ be an $n$-uplet of natural numbers $k_i\geq 1$.
Put
$$
a_I
=a_{k_1}\krest a_{k_2}\krest  ... \krest a_{k_n}.$$
This is a
$(n+1)$-component
link.
\end{defi}

Consider the orientation $\OOOO$ of this link as indicated on the Fig. 1
and let $\mu_1, ... ,\mu_{n+1}$ be the positive
meridians of the components of $a_I$
\wrt~ this orientation. Let $\Gamma_i$ (where $1\leq i \leq n$)
be the integer hyperplane
in $H^1(S^3\sm a_I)$ defined by
$$
\G_i=
\{ \xi~|~ \xi(\m_i)=\xi(\mu_{i+1})\}.
$$
Let $\r(\xi, I )$ be the number of
hyperplanes $\G_i$ \sut~ $\xi\in\G_i$.

\beth\label{t:MN_aI}
Assume that $\xi(\mu_j)\not=0$ for every $j$.
%$ 2\leq j \leq n$.
Then
$$\MM\NN(a_I, \xi)=2\r(\xi).$$
\enth
\Prf
The inequality
$\MM\NN(a_I, \xi)\leq 2\r(\xi)$
follows from the following general property
(the proof is simipar to that of Proposition 6.2 of \cite{PRW}
and will be omitted).
\bele\label{l:conn-sum-constr}
Let $L$ be the
connected  sum of links
$L_1, L_2$ along the components $K_1, K_2$.
Denote by $\mu_i$ the meridians of $K_i$,
here $i=1,2$.
Let $\xi\in H^1(E(L))$, and
assume that
$\xi(\mu_i)\not=0$.
Then
\bq\lb{f:ineq-conn-sum}
\MM\NN(L, \xi) \leq
\MM\NN(L_1,\xi_1)+ \MM\NN(L_2, \xi_2)
\end{equation}
where $\xi_i$ are the restrictions of the class $\xi$ to the
complements of $L_i$.$\qs$
\enle
%is easy since in general
%$\MM\NN(A_1\krest A_2, \xi) \leq
%\MM\NN(A_1,\xi_1)+ \MM\NN(A_2, \xi_2)$
%(where $\xi_i$ are the restrictions of the class $\xi$ to the
%complements of $A_i$).
The proof of the inverse inequality requires
computation of the first Novikov homology of $E(a_I)$
associated to its maximal abelian covering.

Denote by $b_i$ the homology class of the meridian $\mu_i$.
Then $H_1(E(L))$ is freely generated by $b_i$. Let
$\PPPP=\zz[b_1^{\pm 1}, \ldots b_{n+1}^{\pm 1}]$.
Put
$$
Q_m(x,y)=x^m+x^{m-1}y+\ldots +xy^{m-1} + y^m,
\ \ \
S_j=Q_{k_j} (b_j, b_{j+1}).
$$
\bepr\label{p:homol_L}
$$
H_*\big(E(L), \Pxi\bigskip)
\approx
\underset {i=1} {\overset n \bigoplus} \Pxi\Big/ (S_i).
$$
\enpr
\Prf
Denote $\zz[b_i^{\pm 1},  b_{i+1}^{\pm 1}]$ by $\PPPP_i$,
and by $\Pxii$
the Novikov completion of $\PPPP_i$ \wrt~ the restriction $\xi_i$ of $\xi$
to $E(L_i)$.
We have
$$
H_*\big(E(L_i), \Pxii\bigskip)
\approx
\Pxii\Big/ (S_i),
$$
see \rrf{f:qk}.
Consider the maximal abelian covering
$\ove{E(L)} \arrr p E(L)$
and denote by
$\ove{E_i} \arrr {p_i} E(L_i) $ its restriction
to $E(L_i)$.
It is easy to see that
$$
H_*(\ove{E_i} )
\approx
\PPPP\tens{\PPPP_i} H_*(E(L_i)),
\ \
%{\rm and }
\ \
H_*(\ove{E_i}, \Pxi)
\approx
\Pxi\tens{\Pxii} H_*(E(L_i)),
$$
and the proof is over. $\qs$

For a module $N$ over a commutative
 ring $R$
denote by $\mu_R(N)$ the minimal possible number of $R$-generators of $N$.
We will abbreviate this notation to $\mu(N)$ if \noconf.
Observe that for any Morse map
$f:E(L)\to S^1$ belonging to a class $\xi$
we have
$$
m_1(f)\geq \mu\left(H_1(E(L), \Pxi)\right),$$
%In the case when $M=E(L)$ we have
therefore
$$
\MM\NN(E(L),\xi)
\geq
2\mu\left(H_1(E(L), \Pxi\right).
$$
Denote the $\Pxi$-module $\Pxi/(S_j)$ by $W_j$.
The  assertion of Theorem \ref{t:MN_aI} follows from the next Lemma.
\bele\label{l:mu}
$$
\mu\left(\underset {j=1} {\overset n \bigoplus} W_j
\right)
\geq
\r(\xi).
$$
\enle
\noindent
\Prf
We will use in the computation the Novikov ring $\cc((t))$
of Laurent formal power series in one variable with coefficients in $\cc$,
this ring is a field.
Put $a_i=\xi(b_i)\not=0$.
Consider a \ho~
$\phi:\PPPP\to \cc((t))$
defined by $\phi(b_i)=\l_it^{-a_i}$
where $0\not=\l_i\in\cc$.
It is easy to show that $\phi$ is extensible to a \ho~
$\wh\phi:\Pxi\to \cc((t))$.
Choose $\l_j, 1\leq j \leq n+1$ so as to satisfy the following condition:
$$
{\rm if } \ \
a_j=a_{j+1} \ \
{\rm then } \ \
\l_{j+1}/\l_j = \t,
$$
 where $\t\in\cc, \t\not=1$
 is a non-trivial root of unity of degree $k_j$.
% and
%$\t^{k_j}=1$.
Observe that with this choice of $\t$ we have the following property:
\bq\lb{f:theta}
{\rm if } \ \
a_j=a_{j+1} \ \
{\rm then } \ \
\phi(S_j)=0.
\end{equation}

\noindent
The map $\wh\phi$ endows $\cc((t))$ with a structure
of $\Pxi$-module.
%Denote the module
%$\Pxi/(S_j)$ by $W_j$.
The condition \rrf{f:theta} implies
$$
W_j\tens{\wh \phi} \cc((t))\approx \cc((t)) \ \
{\rm if } \ \
a_j=a_{j+1}.
$$
Therefore
$$
\mu(H_1(E(L), \Lxi))
=
\mu\Big(\underset {i=1} {\overset n \bigoplus} W_j
\Big)
\geq
\mu_{\cc((t))}\Big(\underset {i=1} {\overset n \bigoplus} W_j\tens{\wh\phi}
\cc((t))
\Big)
\geq \r(\xi, I).
$$
$\qs$

\section{Links with $\leq 8$ crossings}
 \label{s:6-7-8}

 In this section, we present the results of our computation of the Morse-Novikov numbers for all 2-component links with crossing number less than eight
 (except the link L8a7, for which he have a partial result).
 For a cohomology class $\xi=px^*+qy^*\in H^1(C_L; \mathbb {Z})$, here $x,y$ are the meridional homology classes of $L$.
 we represent its Morse-Novikov number $\mathcal{MN}(L,\xi)$ graphically.
 In these diagrams, each cohomology class $\xi$ corresponds to an intersection point of the unit circle with the line $py^*+qx^*=0$.
 Black points denote the non-regular classes.
 Green points indicate that $2 \leq \mathcal{MN}(L,\xi)\leq 4$,
 Red points indicate that $\mathcal{MN}(L,\xi)=2$,
 while blue points indicate that $\mathcal{MN}(L,\xi)=0$,
 in which case the link is fibered with respect to $\xi$.

For example, the diagram for L5a1 shows that $\xi$ is non-regular when $p=0$ or $q=0$;
otherwise $\mathcal{MN}(L,\xi)=0$.
The diagram for L6a1 shows that $\mathcal{MN}(L,\xi)=0$ for $pq \leqslant 0$ and $\mathcal{MN}(L,\xi)=2$ for $pq>0$.
Finally, the diagram for L7n1 shows that $\mathcal{MN}(L,\xi)=2$ precisely when $3p+q=0$ and $\mathcal{MN}(L,\xi)=0$ otherwise.

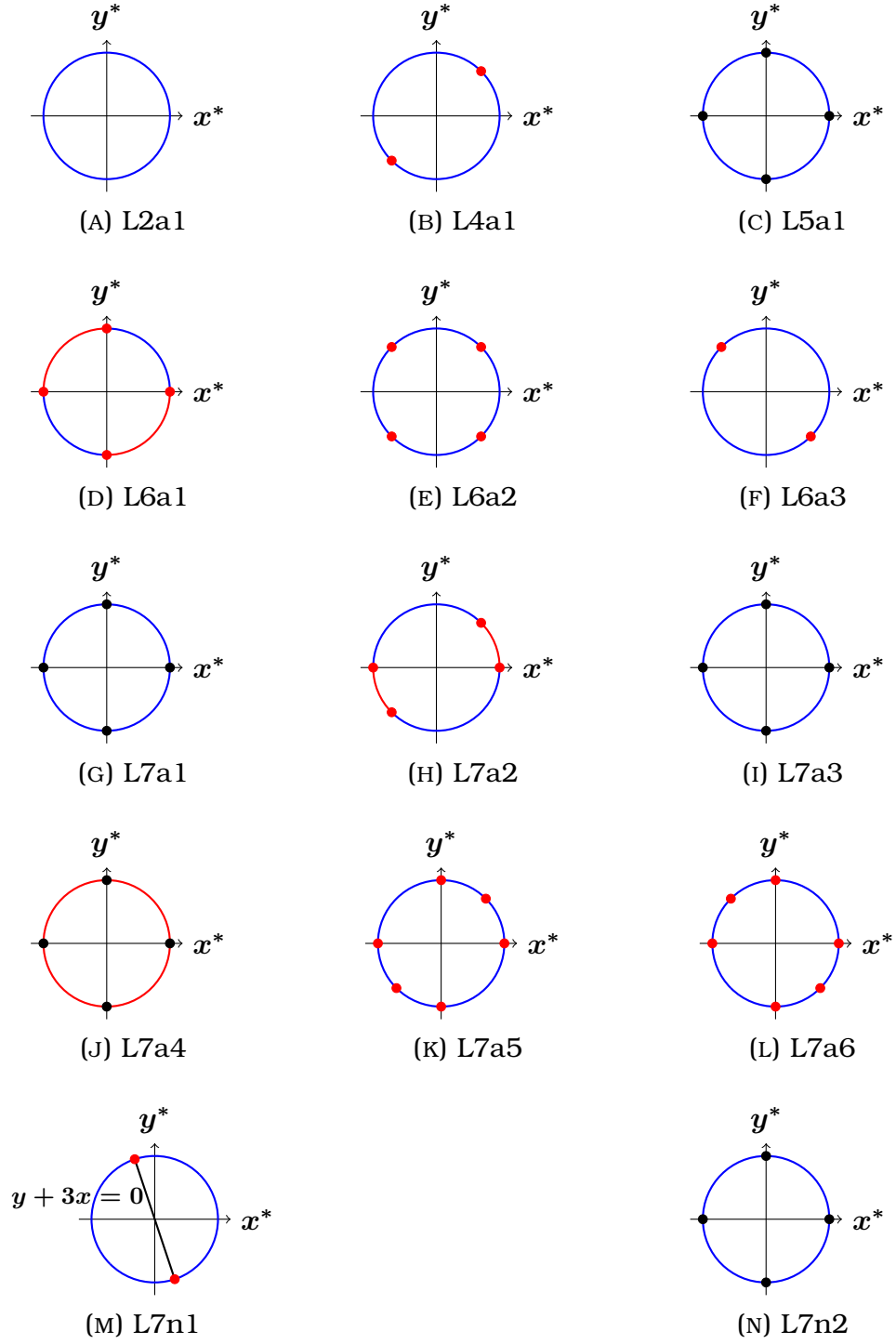
\begin{figure}[p]
      \centering
      \captionsetup[subfigure]{labelformat=parens, labelsep=space, size=small}

      % row 1
      \begin{subfigure}[b]{0.32\textwidth}
          \centering
          \begin{tikzpicture}[scale=0.9]
              % axis
              \draw[->] (-1.2,0) -- (1.2,0) node[right] {$x^*$};
              \draw[->] (0,-1.2) -- (0,1.2) node[above] {$y^*$};

              % circle
              \draw[thick, blue] (1,0) arc (0:360:1);

              % point

          \end{tikzpicture}
          \caption{L2a1}
      \end{subfigure}
      \hfill
      \begin{subfigure}[b]{0.32\textwidth}
          \centering
          \begin{tikzpicture}[scale=0.9]
              % axis
              \draw[->] (-1.2,0) -- (1.2,0) node[right] {$x^*$};
              \draw[->] (0,-1.2) -- (0,1.2) node[above] {$y^*$};

              %

              % circle
              \draw[thick, blue] (1,0) arc (0:360:1);

              % point
              \filldraw[red] (45:1) circle (2pt);
              \filldraw[red] (225:1) circle (2pt);
          \end{tikzpicture}
          \caption{L4a1}
      \end{subfigure}
      \hfill
      \begin{subfigure}[b]{0.32\textwidth}
          \centering
          \begin{tikzpicture}[scale=0.9]
              % axis
              \draw[->] (-1.2,0) -- (1.2,0) node[right] {$x^*$};
              \draw[->] (0,-1.2) -- (0,1.2) node[above] {$y^*$};

              % circle
              \draw[thick, blue] (1,0) arc (0:360:1);

              % point
              \filldraw[black] (0:1) circle (2pt);
              \filldraw[black] (90:1) circle (2pt);
              \filldraw[black] (180:1) circle (2pt);
              \filldraw[black] (270:1) circle (2pt);
          \end{tikzpicture}
          \caption{L5a1}
      \end{subfigure}
      \vspace{0.4cm}

      % row 2
      \begin{subfigure}[b]{0.32\textwidth}
          \centering
          \begin{tikzpicture}[scale=0.9]
              % axis
              \draw[->] (-1.2,0) -- (1.2,0) node[right] {$x^*$};
              \draw[->] (0,-1.2) -- (0,1.2) node[above] {$y^*$};

              % circle
              \draw[thick, blue] (1,0) arc (0:90:1);
              \draw[thick, blue] (-1,0) arc (180:270:1);
              \draw[thick,red] (0,1) arc (90:180:1);
              \draw[thick,red] (0,-1) arc (270:360:1);

              % point
              \filldraw[red] (0:1) circle (2pt);
              \filldraw[red] (90:1) circle (2pt);
              \filldraw[red] (180:1) circle (2pt);
              \filldraw[red] (270:1) circle (2pt);

          \end{tikzpicture}
          \caption{L6a1}
      \end{subfigure}
      \hfill
      \begin{subfigure}[b]{0.32\textwidth}
          \centering
          \begin{tikzpicture}[scale=0.9]
              % axis
              \draw[->] (-1.2,0) -- (1.2,0) node[right] {$x^*$};
              \draw[->] (0,-1.2) -- (0,1.2) node[above] {$y^*$};

              % circle
              \draw[thick, blue] (1,0) arc (0:90:1);
              \draw[thick, blue] (-1,0) arc (180:270:1);
              \draw[thick,blue] (0,1) arc (90:180:1);
              \draw[thick,blue] (0,-1) arc (270:360:1);

              % point
              \filldraw[red] (45:1) circle (2pt);
              \filldraw[red] (135:1) circle (2pt);
              \filldraw[red] (225:1) circle (2pt);
              \filldraw[red] (315:1) circle (2pt);

          \end{tikzpicture}
          \caption{L6a2}
      \end{subfigure}
      \hfill
      \begin{subfigure}[b]{0.32\textwidth}
          \centering
          \begin{tikzpicture}[scale=0.9]
              % axis
              \draw[->] (-1.2,0) -- (1.2,0) node[right] {$x^*$};
              \draw[->] (0,-1.2) -- (0,1.2) node[above] {$y^*$};

              % circle
              \draw[thick, blue] (1,0) arc (0:90:1);
              \draw[thick, blue] (-1,0) arc (180:270:1);
              \draw[thick,blue] (0,1) arc (90:180:1);
              \draw[thick,blue] (0,-1) arc (270:360:1);

              % point

              \filldraw[red] (135:1) circle (2pt);

              \filldraw[red] (315:1) circle (2pt);

          \end{tikzpicture}
          \caption{L6a3}
      \end{subfigure}
      \vspace{0.4cm}

      % row 3
      \begin{subfigure}[b]{0.32\textwidth}
          \centering
          \begin{tikzpicture}[scale=0.9]
              % axis
              \draw[->] (-1.2,0) -- (1.2,0) node[right] {$x^*$};
              \draw[->] (0,-1.2) -- (0,1.2) node[above] {$y^*$};

              % circle
              \draw[thick, blue] (1,0) arc (0:90:1);
              \draw[thick, blue] (-1,0) arc (180:270:1);
              \draw[thick,blue] (0,1) arc (90:180:1);
              \draw[thick,blue] (0,-1) arc (270:360:1);

              % point
              \filldraw[black] (0:1) circle (2pt);
              \filldraw[black] (90:1) circle (2pt);
              \filldraw[black] (180:1) circle (2pt);
              \filldraw[black] (270:1) circle (2pt);

          \end{tikzpicture}
          \caption{L7a1}
      \end{subfigure}
      \hfill
      \begin{subfigure}[b]{0.32\textwidth}
          \centering
          \begin{tikzpicture}[scale=0.9]
              % axis
              \draw[->] (-1.2,0) -- (1.2,0) node[right] {$x^*$};
              \draw[->] (0,-1.2) -- (0,1.2) node[above] {$y^*$};

              % circle
              \draw[thick, red] (1,0) arc (0:45:1);
              \draw[thick, red] (-1,0) arc (180:225:1);
              \draw[thick,blue] (45:1) arc (45:180:1);
              \draw[thick,blue] (225:1) arc (225:360:1);

              % point
              \filldraw[red] (45:1) circle (2pt);
              \filldraw[red] (0:1) circle (2pt);
              \filldraw[red] (225:1) circle (2pt);
              \filldraw[red] (180:1) circle (2pt);

          \end{tikzpicture}
          \caption{L7a2}
      \end{subfigure}
      \hfill
      \begin{subfigure}[b]{0.32\textwidth}
          \centering
          \begin{tikzpicture}[scale=0.9]
              % axis
              \draw[->] (-1.2,0) -- (1.2,0) node[right] {$x^*$};
              \draw[->] (0,-1.2) -- (0,1.2) node[above] {$y^*$};

              % circle
              \draw[thick, blue] (1,0) arc (0:90:1);
              \draw[thick, blue] (-1,0) arc (180:270:1);
              \draw[thick,blue] (0,1) arc (90:180:1);
              \draw[thick,blue] (0,-1) arc (270:360:1);

              % point
              \filldraw[black] (0:1) circle (2pt);
              \filldraw[black] (90:1) circle (2pt);
              \filldraw[black] (180:1) circle (2pt);
              \filldraw[black] (270:1) circle (2pt);

          \end{tikzpicture}
          \caption{L7a3}
      \end{subfigure}
      \vspace{0.4cm}

      % row 4
      \begin{subfigure}[b]{0.32\textwidth}
          \centering
          \begin{tikzpicture}[scale=0.9]
              % axis
              \draw[->] (-1.2,0) -- (1.2,0) node[right] {$x^*$};
              \draw[->] (0,-1.2) -- (0,1.2) node[above] {$y^*$};

              % circle
              \draw[thick, red] (1,0) arc (0:90:1);
              \draw[thick, red] (-1,0) arc (180:270:1);
              \draw[thick,red] (0,1) arc (90:180:1);
              \draw[thick,red] (0,-1) arc (270:360:1);

              % point
              \filldraw[black] (0:1) circle (2pt);
              \filldraw[black] (90:1) circle (2pt);
              \filldraw[black] (180:1) circle (2pt);
              \filldraw[black] (270:1) circle (2pt);

          \end{tikzpicture}
          \caption{L7a4}
      \end{subfigure}
      \hfill
      \begin{subfigure}[b]{0.32\textwidth}
          \centering
          \begin{tikzpicture}[scale=0.9]
              % axis
              \draw[->] (-1.2,0) -- (1.2,0) node[right] {$x^*$};
              \draw[->] (0,-1.2) -- (0,1.2) node[above] {$y^*$};

              % circle
              \draw[thick, blue] (1,0) arc (0:45:1);
              \draw[thick,blue] (-1,0) arc (180:225:1);
              \draw[thick,blue] (45:1) arc (45:180:1);
              \draw[thick,blue] (225:1) arc (225:360:1);

              % point
              \filldraw[red] (90:1) circle (2pt);
              \filldraw[red] (270:1) circle (2pt);
              \filldraw[red] (45:1) circle (2pt);
              \filldraw[red] (0:1) circle (2pt);
              \filldraw[red] (225:1) circle (2pt);
              \filldraw[red] (180:1) circle (2pt);

          \end{tikzpicture}
          \caption{L7a5}
      \end{subfigure}
      \hfill
      \begin{subfigure}[b]{0.32\textwidth}
          \centering
          \begin{tikzpicture}[scale=0.9]
              % axis
              \draw[->] (-1.2,0) -- (1.2,0) node[right] {$x^*$};
              \draw[->] (0,-1.2) -- (0,1.2) node[above] {$y^*$};

              % circle
              \draw[thick, blue] (1,0) arc (0:45:1);
              \draw[thick,blue] (-1,0) arc (180:225:1);
              \draw[thick,blue] (45:1) arc (45:180:1);
              \draw[thick,blue] (225:1) arc (225:360:1);

              % point
              \filldraw[red] (90:1) circle (2pt);
              \filldraw[red] (270:1) circle (2pt);
              \filldraw[red] (135:1) circle (2pt);
              \filldraw[red] (0:1) circle (2pt);
              \filldraw[red] (315:1) circle (2pt);
              \filldraw[red] (180:1) circle (2pt);

          \end{tikzpicture}
          \caption{L7a6}
      \end{subfigure}\vspace{0.4cm}

      % row 5
      \makebox[\textwidth][c]{
	  \begin{subfigure}[b]{0.32\textwidth}

    	\centering
   		 \begin{tikzpicture}[scale=0.9]
        % axis
        \draw[->] (-1.2,0) -- (1.2,0) node[right] {$x^*$};
        \draw[->] (0,-1.2) -- (0,1.2) node[above] {$y^*$};

        % circle arcs
        \draw[thick, blue] (1,0) arc (0:45:1);
        \draw[thick,blue] (-1,0) arc (180:225:1);
        \draw[thick,blue] (45:1) arc (45:180:1);
        \draw[thick,blue] (225:1) arc (225:360:1);

        \draw[thick, black!70!black]
            ({(1/sqrt(10))}, {(-3/sqrt(10))}) -- ({(-1/sqrt(10))}, {(3/sqrt(10))})
            node[midway, above left] {\footnotesize $y+3x=0$};

        \filldraw[red] ({(1/sqrt(10))}, {(-3/sqrt(10))}) circle (2pt);
        \filldraw[red] ({(-1/sqrt(10))}, {(3/sqrt(10))}) circle (2pt);

    	\end{tikzpicture}
   		\caption{L7n1}
   	\end{subfigure}
      \hfill
      \begin{subfigure}[b]{0.32\textwidth}
          \centering
          \begin{tikzpicture}[scale=0.9]
              % axis
              \draw[->] (-1.2,0) -- (1.2,0) node[right] {$x^*$};
              \draw[->] (0,-1.2) -- (0,1.2) node[above] {$y^*$};

              % circle
              \draw[thick, blue] (1,0) arc (0:90:1);
              \draw[thick, blue] (-1,0) arc (180:270:1);
              \draw[thick,blue] (0,1) arc (90:180:1);
              \draw[thick,blue] (0,-1) arc (270:360:1);

              % point
              \filldraw[black] (0:1) circle (2pt);
              \filldraw[black] (90:1) circle (2pt);
              \filldraw[black] (180:1) circle (2pt);
              \filldraw[black] (270:1) circle (2pt);

          \end{tikzpicture}
          \caption{L7n2}
      \end{subfigure}

      }

      \caption{Morse-Novikov number with respect to $\xi$}
      \label{fig:links-row3}

\end{figure}

\begin{figure}[p]
\centering
\captionsetup[subfigure]{labelformat=parens, labelsep=space, size=small}

% row 1
\begin{subfigure}[b]{0.32\textwidth}
    \centering
    \begin{tikzpicture}[scale=0.9]
        % axis
        \draw[->] (-1.2,0) -- (1.2,0) node[right] {$x^*$};
        \draw[->] (0,-1.2) -- (0,1.2) node[above] {$y^*$};

        % circle
        \draw[thick, blue] (1,0) arc (0:90:1);
        \draw[thick, blue] (-1,0) arc (180:270:1);
        \draw[thick,blue] (0,1) arc (90:180:1);
        \draw[thick,blue] (0,-1) arc (270:360:1);

        % point
        \filldraw[black] (0:1) circle (2pt);
        \filldraw[black] (90:1) circle (2pt);
        \filldraw[black] (180:1) circle (2pt);
        \filldraw[black] (270:1) circle (2pt);

    \end{tikzpicture}
    \caption{L8a1}
\end{subfigure}
\hfill
\begin{subfigure}[b]{0.32\textwidth}
    \centering
    \begin{tikzpicture}[scale=0.9]
        % axis
        \draw[->] (-1.2,0) -- (1.2,0) node[right] {$x^*$};
        \draw[->] (0,-1.2) -- (0,1.2) node[above] {$y^*$};

        % circle
        \draw[thick, blue] (1,0) arc (0:90:1);
        \draw[thick, blue] (-1,0) arc (180:270:1);
        \draw[thick,blue] (0,1) arc (90:180:1);
        \draw[thick,blue] (0,-1) arc (270:360:1);

        % point
        \filldraw[black] (0:1) circle (2pt);
        \filldraw[black] (90:1) circle (2pt);
        \filldraw[black] (180:1) circle (2pt);
        \filldraw[black] (270:1) circle (2pt);

    \end{tikzpicture}
    \caption{L8a2}
\end{subfigure}
\hfill
\begin{subfigure}[b]{0.32\textwidth}
    \centering
    \begin{tikzpicture}[scale=0.9]
        % axis
        \draw[->] (-1.2,0) -- (1.2,0) node[right] {$x^*$};
        \draw[->] (0,-1.2) -- (0,1.2) node[above] {$y^*$};

        % circle
        \draw[thick, red] (1,0) arc (0:45:1);
        \draw[thick,red] (-1,0) arc (180:225:1);
        \draw[thick,blue] (45:1) arc (45:180:1);
        \draw[thick,blue] (225:1) arc (225:360:1);

        % point

        \filldraw[red] (45:1) circle (2pt);
        \filldraw[red] (0:1) circle (2pt);
        \filldraw[red] (225:1) circle (2pt);
        \filldraw[red] (180:1) circle (2pt);

    \end{tikzpicture}
    \caption{L8a3}
\end{subfigure}\vspace{0.4cm}

% row 2
\begin{subfigure}[b]{0.32\textwidth}
    \centering
    \begin{tikzpicture}[scale=0.9]
        % axis
        \draw[->] (-1.2,0) -- (1.2,0) node[right] {$x^*$};
        \draw[->] (0,-1.2) -- (0,1.2) node[above] {$y^*$};

        % circle
        \draw[thick, blue] (1,0) arc (0:90:1);
        \draw[thick, blue] (-1,0) arc (180:270:1);
        \draw[thick,blue] (0,1) arc (90:180:1);
        \draw[thick,blue] (0,-1) arc (270:360:1);

        % point
        \filldraw[black] (0:1) circle (2pt);
        \filldraw[black] (90:1) circle (2pt);
        \filldraw[black] (180:1) circle (2pt);
        \filldraw[black] (270:1) circle (2pt);

    \end{tikzpicture}
    \caption{L8a4}
\end{subfigure}
\hfill
\begin{subfigure}[b]{0.32\textwidth}
    \centering
    \begin{tikzpicture}[scale=0.9]
        % axis
        \draw[->] (-1.2,0) -- (1.2,0) node[right] {$x^*$};
        \draw[->] (0,-1.2) -- (0,1.2) node[above] {$y^*$};

        % circle
        \draw[thick, blue] (1,0) arc (0:90:1);
        \draw[thick, blue] (-1,0) arc (180:270:1);
        \draw[thick,red] (0,1) arc (90:180:1);
        \draw[thick,red] (0,-1) arc (270:360:1);

        % point
        \filldraw[red] (0:1) circle (2pt);
        \filldraw[red] (90:1) circle (2pt);
        \filldraw[red] (180:1) circle (2pt);
        \filldraw[red] (270:1) circle (2pt);

    \end{tikzpicture}
    \caption{L8a5}
\end{subfigure}
\hfill
\begin{subfigure}[b]{0.32\textwidth}
    \centering
    \begin{tikzpicture}[scale=0.9]
        % axis
        \draw[->] (-1.2,0) -- (1.2,0) node[right] {$x^*$};
        \draw[->] (0,-1.2) -- (0,1.2) node[above] {$y^*$};

        % circle
        \draw[thick, red] (1,0) arc (0:45:1);
        \draw[thick,red] (-1,0) arc (180:225:1);
        \draw[thick,red] (45:1) arc (45:180:1);
        \draw[thick,red] (225:1) arc (225:360:1);

        % point

    \end{tikzpicture}
    \caption{L8a6}
\end{subfigure}\vspace{0.4cm}

% row 3
\begin{subfigure}[b]{0.32\textwidth}
    \centering
    \begin{tikzpicture}[scale=0.9]
        % axis
        \draw[->] (-1.2,0) -- (1.2,0) node[right] {$x^*$};
        \draw[->] (0,-1.2) -- (0,1.2) node[above] {$y^*$};

        % circle
        \draw[thick,green] (1,0) arc (0:45:1);
        \draw[thick,green] (-1,0) arc (180:225:1);
        \draw[thick,blue] (45:1) arc (45:180:1);
        \draw[thick,blue] (225:1) arc (225:360:1);

        % point

        \filldraw[red] (45:1) circle (2pt);
        \filldraw[green] (0:1) circle (2pt);
        \filldraw[red] (225:1) circle (2pt);
        \filldraw[green] (180:1) circle (2pt);

    \end{tikzpicture}
    \caption{L8a7}
    \end{subfigure}
\hfill
\begin{subfigure}[b]{0.32\textwidth}
    \centering
    \begin{tikzpicture}[scale=0.9]
        % axis
        \draw[->] (-1.2,0) -- (1.2,0) node[right] {$x^*$};
        \draw[->] (0,-1.2) -- (0,1.2) node[above] {$y^*$};

        % circle
        \draw[thick, blue] (1,0) arc (0:90:1);
        \draw[thick, blue] (-1,0) arc (180:270:1);
        \draw[thick,blue] (0,1) arc (90:180:1);
        \draw[thick,blue] (0,-1) arc (270:360:1);

        % point
        \filldraw[red] (0:1) circle (2pt);
        \filldraw[red] (90:1) circle (2pt);
        \filldraw[red] (180:1) circle (2pt);
        \filldraw[red] (270:1) circle (2pt);

    \end{tikzpicture}
    \caption{L8a8}
\end{subfigure}
\hfill
\begin{subfigure}[b]{0.32\textwidth}
    \centering
    \begin{tikzpicture}[scale=0.9]
        % axis
        \draw[->] (-1.2,0) -- (1.2,0) node[right] {$x^*$};
        \draw[->] (0,-1.2) -- (0,1.2) node[above] {$y^*$};

        % circle
        \draw[thick, blue] (1,0) arc (0:90:1);
        \draw[thick, blue] (-1,0) arc (180:270:1);
        \draw[thick,blue] (0,1) arc (90:180:1);
        \draw[thick,blue] (0,-1) arc (270:360:1);

        % point
        \filldraw[red] (0:1) circle (2pt);
        \filldraw[red] (90:1) circle (2pt);
        \filldraw[red] (180:1) circle (2pt);
        \filldraw[red] (270:1) circle (2pt);

    \end{tikzpicture}
    \caption{L8a9}

\end{subfigure}\vspace{0.4cm}

% row 4
\begin{subfigure}[b]{0.32\textwidth}
    \centering
    \begin{tikzpicture}[scale=0.9]
        % axis
        \draw[->] (-1.2,0) -- (1.2,0) node[right] {$x^*$};
        \draw[->] (0,-1.2) -- (0,1.2) node[above] {$y^*$};

        % circle
        \draw[thick, red] (1,0) arc (0:90:1);
        \draw[thick, red] (-1,0) arc (180:270:1);
        \draw[thick,blue] (0,1) arc (90:180:1);
        \draw[thick,blue] (0,-1) arc (270:360:1);

        % point
        \filldraw[red] (0:1) circle (2pt);
        \filldraw[red] (90:1) circle (2pt);
        \filldraw[red] (180:1) circle (2pt);
        \filldraw[red] (270:1) circle (2pt);

    \end{tikzpicture}
    \caption{L8a10}
    \end{subfigure}
\hfill
\begin{subfigure}[b]{0.32\textwidth}
    \centering
    \begin{tikzpicture}[scale=0.9]
        % axis
        \draw[->] (-1.2,0) -- (1.2,0) node[right] {$x^*$};
        \draw[->] (0,-1.2) -- (0,1.2) node[above] {$y^*$};

        % circle
        \draw[thick, red] (1,0) arc (0:90:1);
        \draw[thick, red] (-1,0) arc (180:270:1);
        \draw[thick,blue] (0,1) arc (90:180:1);
        \draw[thick,blue] (0,-1) arc (270:360:1);

        % point
        \filldraw[red] (0:1) circle (2pt);
        \filldraw[red] (90:1) circle (2pt);
        \filldraw[red] (180:1) circle (2pt);
        \filldraw[red] (270:1) circle (2pt);
        \filldraw[red] (135:1) circle (2pt);
        \filldraw[red] (315:1) circle (2pt);

    \end{tikzpicture}
    \caption{L8a11}
\end{subfigure}
\hfill
\begin{subfigure}[b]{0.32\textwidth}
    \centering
    \begin{tikzpicture}[scale=0.9]
        % axis
        \draw[->] (-1.2,0) -- (1.2,0) node[right] {$x^*$};
        \draw[->] (0,-1.2) -- (0,1.2) node[above] {$y^*$};

        % circle
        \draw[thick, blue] (1,0) arc (0:90:1);
        \draw[thick, blue] (-1,0) arc (180:270:1);
        \draw[thick,blue] (0,1) arc (90:180:1);
        \draw[thick,blue] (0,-1) arc (270:360:1);

        % point
        \filldraw[red] (45:1) circle (2pt);
        \filldraw[red] (135:1) circle (2pt);
        \filldraw[red] (225:1) circle (2pt);
        \filldraw[red] (315:1) circle (2pt);

    \end{tikzpicture}
    \caption{L8a12}

\end{subfigure}\vspace{0.4cm}

%row5
\makebox[\textwidth][c]{
\begin{subfigure}[b]{0.32\textwidth}
    \centering
    \begin{tikzpicture}[scale=0.9]
        % axis
        \draw[->] (-1.2,0) -- (1.2,0) node[right] {$x^*$};
        \draw[->] (0,-1.2) -- (0,1.2) node[above] {$y^*$};

        % circle
        \draw[thick, blue] (1,0) arc (0:90:1);
        \draw[thick, blue] (-1,0) arc (180:270:1);
        \draw[thick,blue] (0,1) arc (90:180:1);
        \draw[thick,blue] (0,-1) arc (270:360:1);

        % point
        \filldraw[red] (45:1) circle (2pt);
        \filldraw[red] (135:1) circle (2pt);
        \filldraw[red] (225:1) circle (2pt);
        \filldraw[red] (315:1) circle (2pt);

    \end{tikzpicture}
    \caption{L8a13}
    \end{subfigure}
\hfill
\begin{subfigure}[b]{0.32\textwidth}
    \centering
    \begin{tikzpicture}[scale=0.9]
        % axis
        \draw[->] (-1.2,0) -- (1.2,0) node[right] {$x^*$};
        \draw[->] (0,-1.2) -- (0,1.2) node[above] {$y^*$};

        % circle
        \draw[thick, blue] (1,0) arc (0:90:1);
        \draw[thick, blue] (-1,0) arc (180:270:1);
        \draw[thick,blue] (0,1) arc (90:180:1);
        \draw[thick,blue] (0,-1) arc (270:360:1);

        % point
        %\filldraw[red] (0:1) circle (2pt);
        %\filldraw[red] (90:1) circle (2pt);
        %\filldraw[red] (180:1) circle (2pt);
        %\filldraw[red] (270:1) circle (2pt);
        \filldraw[red] (135:1) circle (2pt);
        \filldraw[red] (315:1) circle (2pt);

    \end{tikzpicture}
    \caption{L8a14}
\end{subfigure}

}
%row6
\makebox[\textwidth][c]{
\begin{subfigure}[b]{0.32\textwidth}
    \centering
 \begin{tikzpicture}[scale=0.9]
        % axis
        \draw[->] (-1.2,0) -- (1.2,0) node[right] {$x^*$};
        \draw[->] (0,-1.2) -- (0,1.2) node[above] {$y^*$};

        \draw[thick, blue] (0:1) arc (0:153.435:1);

        \draw[thick, blue] (180:1) arc (180:333.435:1);

        \draw[thick, red] (153.435:1) arc (153.435:180:1);
        \draw[thick, red] (333.435:1) arc (333.435:360:1);

        \draw[thick, black]
            ({2/sqrt(5)}, {-1/sqrt(5)}) -- ({-2/sqrt(5)}, {1/sqrt(5)})
            node[midway, above left] {\footnotesize $x+2y=0$};

        \filldraw[red] (1,0) circle (2pt);   % (1,0)
        \filldraw[red] (-1,0) circle (2pt);  % (-1,0)

        \filldraw[red] ({-2/sqrt(5)}, {1/sqrt(5)}) circle (2pt); % P1
        \filldraw[red] ({2/sqrt(5)}, {-1/sqrt(5)}) circle (2pt);  % P2

    \end{tikzpicture}
    \caption{L8n1}

\end{subfigure}\vspace{0.4cm}
\hfill
\begin{subfigure}[b]{0.32\textwidth}
    \centering
    \begin{tikzpicture}[scale=0.9]
        % axis
        \draw[->] (-1.2,0) -- (1.2,0) node[right] {$x^*$};
        \draw[->] (0,-1.2) -- (0,1.2) node[above] {$y^*$};

        % circle
        \draw[thick, blue] (1,0) arc (0:90:1);
        \draw[thick, blue] (-1,0) arc (180:270:1);
        \draw[thick,blue] (0,1) arc (90:180:1);
        \draw[thick,blue] (0,-1) arc (270:360:1);

        % point
        \filldraw[black] (0:1) circle (2pt);
        \filldraw[black] (90:1) circle (2pt);
        \filldraw[black] (180:1) circle (2pt);
        \filldraw[black] (270:1) circle (2pt);

    \end{tikzpicture}
    \caption{L8n2}
\end{subfigure}

}

\caption{Morse-Novikov number with respect to $\xi$}
%\label{fig:links-row3}
\end{figure}
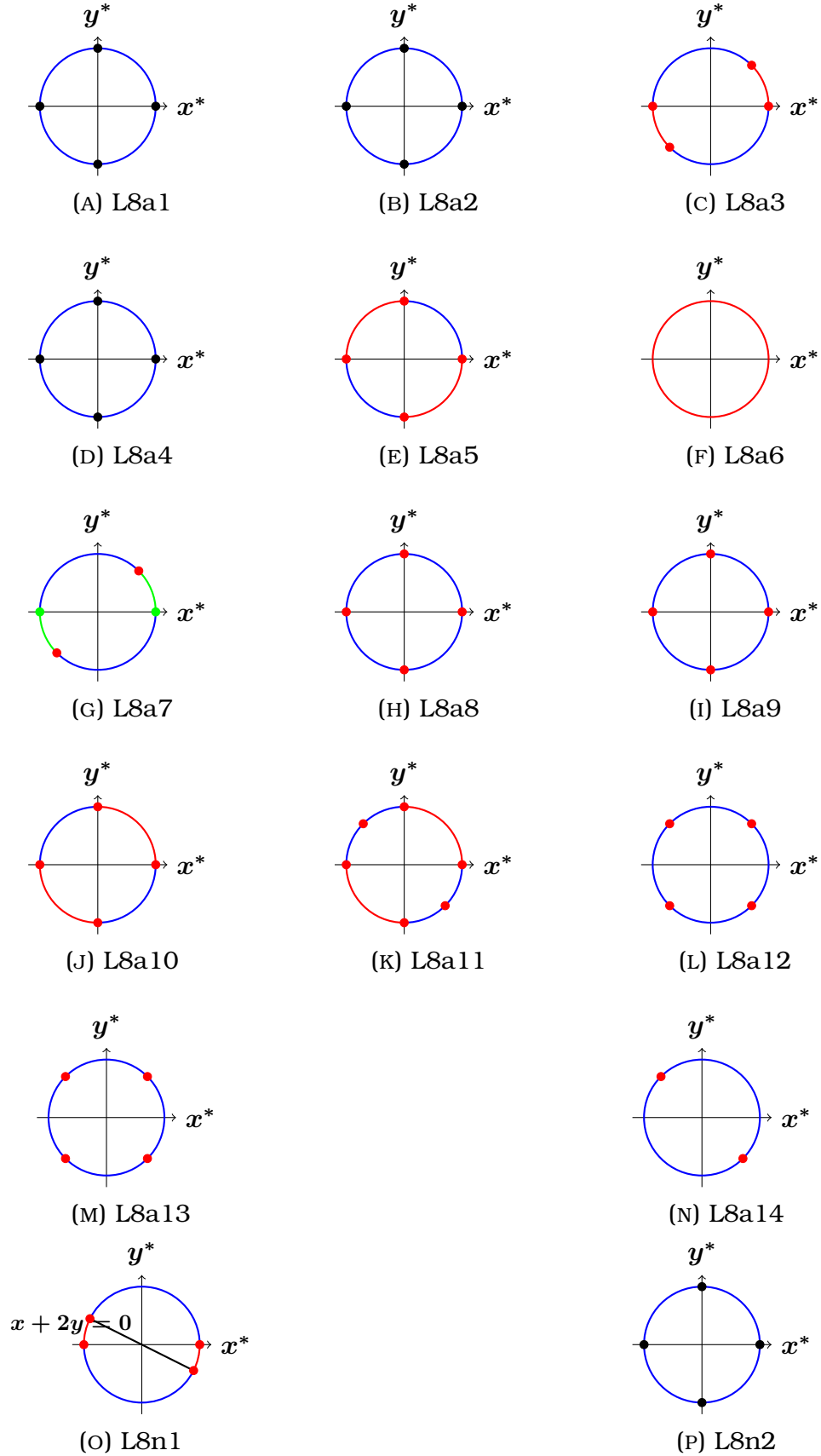

\newpage
\section{The case of closed 1-forms}
\label{s:irrat}

Morse-Novikov theory has a natural generalization to the case
of Morse forms, that is, closed 1-forms with non-degenerate zeros.
The main aim of this section are  Theorems \ref{t:m-vers-mn-irrat} and \ref{t:tunn-ineq}
which generalize Theorems \ref{t:m-vers-mn} and \ref{t:tunnel-forms}
%In this section we establish a generalization of Theorem \ref{??}
to the  case of Morse forms. The generalization of other results
of the previous sections is beyond the scope of the present paper.

We need some more definitions.
\bede \label{d:mores-forms}
Let $M$ be a compact manifold with boundary $\pr M$.
A closed 1-form $\o$ on $M$ is called {\it a Morse form}
if all its zeros are non-degenerate and belong to $M\sm \pr M$,
and the restriction of $\o$ to $\pr M$ is a non-singular 1-form
on $\pr M$.
\end{defi}

For any Morse map $f:M\to S^1$ the differential $df$ is a Morse form.
Such forms will be caled {\it integral }.
A Morse form $\o$ is called {\it rational}, if $\o=C\eta$
where $C\in\qq$ and $\eta$ is an integral form.
It is easy to see that a Morse form is integral, resp. rational
\ifff~ its deRham cohomology  class
$[\omega ]$
is in the image
of
$\hmz \to \hmr$, resp. $\hmq\to\hmr$.

For a Morse form $\o$ we denote by $m_i(\o)$ the number of zeros of $\o$ of index $i$, and by $m(\o)$ the number of all its zeros.

%\bede\label{d:reg-forms}
A class $\xi\in\hmr$ is called {\it regular} if its restriction to $\pr M$ is represented by a closed 1-form without zeros.
A standard Morse-theoretical argument
implies that for any regular class $\xi$ there is a Morse map
$f:M\to S^1$ belonging to $\xi$.
\bede\label{d:mn-forms}
Let $\xi\in\hmr$ be a regular cohomology class. {\it The Morse-Novikov
number  $\MM\NN(M, \xi)$ } of  $\xi$ is the minimal possible number $m(\o)$ where $[\o]=\xi$.
\end{defi}

The next approximation lemma will be used in the proof of
Theorem \ref{t:m-vers-mn-irrat}.

\bele\label{l:approx}
Let $\o\in\O^1(M)$ be a  Morse form,
belonging to a regular class $\xi\in\hmr\approx \rr^n$.
Then in every \nei~ of $\xi\in\hmr$
there is a regular rational cohomology class $\eta\in\hmq$, such that
some integer multiple $N\eta$ is representable by a Morse
map
$\phi:M\to S^1$
having the same Morse numbers as $\o$.
\enle
\Prf
Choose a basis $\{a_1, \ldots , a_n\}$
in the group $H^1(M,\zz)\approx \zz^n$ and represent
each $a_i$ by a closed 1-form $\a_i$ having the following properties:

\been\item
The restriction of $\a_i$ to $\pr M$ is non-singular.
\item
The form $\a_i$ equals $0$ in some \nei~ $U$ of $S(\o)$.
\enen
For a vector $x=(x_1, \ldots , x_n)\in\rr^n$
consider a 1-form
$$
\l(x)
=
\o+\sum_i x_i\a_i.
$$
For every $x$ the 1-form $\l(x)$ equals $\o$
in $U$.
If $x$ is sufficiently small, then $\l(x)$ has no zeros in $M\sm U$, and $\l(x)~|~\pr M$ has no zeros either. Therefore $\l(x)$ is a Morse form and
we have $m_i(\l(x))=m_i(\o)$ for every $x$
sufficiently small. Write $\xi=\sum_i\xi_ia_i$. Choosing $x$
in such a way that $x_i-\xi_i$ be rational, we will have $[\l(x)]\in\hmq$,
so that for some $N\in\nn$ we have
$$
[N\l(x)]\in \hmz
$$
and $N\l(x)$ is the differential of a Morse map to a circle,
naving the same Morse numbers as $\o$.
$\qs$

\beth\label{t:m-vers-mn-irrat}
Let
$M$ be a connected compact manifold with $\pr M\not=\ems$,
and $\xi$ be a regular class in $H^1(M,\rr)$.
Then
$$
\MM\NN(M,\xi)\leq\MM(M; \pr M)-2.
$$
\enth
\Prf
Let $\o$ be a regular Morse form belonging to $\xi$.
The previous lemma implies existence of a regular Morse map $\phi:M\to S^1$
having the same Morse numbers as $\o$.
Apply Theorem \ref{t:m-vers-mn} to complete  the proof. $\qs$

Consider now the case when $M=E(L)$ where $L$ is a link in $S^3$.
A class $\xi\in\hmr$ is regular \ifff~ its image in $H^1(\pr E(K_i),\rr)$
is non-zero for every connected component $K_i$ of the link $L$.

\beth\label{t:tunnel-forms}
Let $\xi\in H^1(\pr E(K_i),\rr)$
be a regular cohomology class. Then
$$\MM\NN(E(L),\xi)
\leq 2t(L).
$$
\enth
\Prf
Recall that $\MM(E(L); \pr E(L))=2t(L)+2.$
Apply the previous theorem and the proof is over.$\qs$

 \newpage
 \section{Conjectures}
 \label{s:conj}

 \subsection{Morse-Novikov  numbers and Thurston cones}

 Let $H$ be a finite dimensional vector space over $\qq$,
 and $\l:H\to\qq$ a non-zero \ho.
 The set
 $\{x\in H ~|~ \l(x)>0\}$
 will be called {\it an
  open half-space}.
 A finite intersection of open half-spaces
 will be called {\it basic convex cone}.
 \bede\label{d:gamma-cone}
A set $A\sbs H$ is called
{\it an A-cone} if
there is a finite family  of hyperplanes $S_k$ in $H$
\sut~$A\sm \cup_k S_k$ is a finite union of basic convex cones.
\end{defi}

For $H=H^1(M^3,\qq)$ the cone on each open face of the Thurston polyhedron
is a  basic convex cone. We will call these cones {\it Thurston cones}.
Thurston's theorem implies in particular that the set of all fibered
classes $\xi\in H^1(M^3,\qq)$ is an $I$-cone
\footnote{By certain abuse of notation we say that
a rational cohomology class $\xi\in H^1(M,\qq)$
is fibered if for some $N\in\nn$ the class
$N\xi$ is an integer fibered cohomology class.}

\begin{conj}\label{conj:G}
 Let $L$ be a link in $S^3$. Then for every
 natural $k$ the set
 $$
 \G_k=\{\xi\in H^1(E(L),\qq) ~|~ \xi \ \ {\rm is\ regular\ and  } \ \ \MM\NN(E(L)), \xi)=k\}
 $$
 is an $A$-cone.
 \end{conj}

 One motivation for this conjecture comes from the following result concerning
 the Novikov homology. Denote by $NH(M,\xi)$ the number

\bq\label{f:mn-eq22}
\sum_{i} \left( \wh b_i (M,\xi) +\wh
 q_i (M,\xi) +
  \wh q_{i-1} (M,\xi) \right),
 \end{equation}
then we have $NH(M,\xi)\leq \MNxi$.

\beth\cite{P-Laurent}
 For a finite CW-complex $M$
 the set
 $$
 \D_k=\{\xi\in H^1(M,\qq) ~|~ \xi\not=0  \ \ and   \ \ NH(M, \xi)=k\}
 $$
 is an $A$-cone.
 \end{theo}

Assuming that the conjecture \ref{conj:G}
holds, we have the following natural question:

{\bf Question.} What is the relation between the cones $\G_k$ and Thurston cones?

\subsection{Morse-Novikov number of connected sum of links}

The Morse-Novikov number of knots is additive:
$$
\MM\NN(K_1\krest K_2)
=
\MM\NN(K_1)+\MM\NN(K_2).
$$
This was prover in \cite{ManjarrezGutierrez} for the case of small knots
and in \cite{Baker} in full generality.

\begin{conj}
 Let $L_1, L_2$ be links, put $L=L_1\krest L_2$ and let $\xi\in H^1(E(L))$
 be a regular class such that its restrictions $\xi_1, \xi_2$ to
  $E(L_1)$ and $E(L_2)$
 are regular. Then
 $$
\MM\NN(L_1\krest L_2,\xi)
=
\MM\NN(L_1,\xi_1)+\MM\NN(L_2, \xi_2).
$$
 \end{conj}

\newpage
\bibliography{mybib}
\bibliographystyle{plain}

\end{document}

In Section \ref{s:mn-manif} we introduce a notion of a Morse-Novikov number
$\MM\NN(M,\xi)$ where $M$ is a compact manifold with boundary and
$\xi$ is a 1-dimensional cohomology class.
We prove that
\bq\label{f:ineqq2}
\MM\NN(M,\xi)
\leq \MM(M)-2
\end{equation}
where $\MM(M)$ is the Morse number of the manifold, that is, the minimal possible number of critivcal points of a Morse function on $M$.
This inequality generalizes the inequality
\rrf{f:citing-tunn} and its analog for 2-knots and sufrace-links,
used in \cite{EP}.
In Subsection \ref{su:mn-links} we prove that for a link $L$ in the 3-sphere we have
$$
\MM\NN(E(L)),\xi)
\leq t(L).
$$

For the case of 3-manifolds the invariant $\MM\NN(M,\xi)$
is presumably equal to the invariant $h(M,\g,\xi)$ introduced
by K. Baker and F. Majarrez-Gutierrez
in \cite{BakerManjarrezGutierrez}, for some special choice of the suture $\g$.

Some more notions  from  \mntt~ (Section
 \ref{s:nov-rings}) are necessary to deal with the case of
\tbtc~ (Section \ref{s:two-br}). We prove that for such a link
$L$
a class $\xi\in  H^1(E(L))$ is fibered
\ifff~ the 2-variable Alexander polynomial of $L$ is $\xi$-monic.
Together with tunnel estimate from Section  \ref{s:tunnel} this allows us to
determine completely the
\mnn s of \tbtc s.

Section \ref{s:conn-sum}
is about connected sums of links.
Let
$L=L_1\krest L_2$ be a connected sum of links $L_1$ and $L_2$,
and $\xi\in H^1(E(L))$.
Under a mild restriction on $\xi$ we prove that
$\xi$ is fibered \ifff~ both its restrictions to $E(L_i), i=1,2$ are fibered.
This generalizes a theorem of D. Gabai \cite{Gabai}
which says that a connected sum of links is fibered \ifff~
each of the links is fibered. Our proof is based on a short computation
with  Novikov homology and does not use the sutured manifold theory.

In Section \ref{s:a_k} we consider one particular family
of links. The 2-component link $a_k$ is defined as the boundary
of a twisted annulus embedded in $S^3$. We compute its \mnn, and
also the \mnn s of connected sums of such links.

In Section
 \ref{s:6-7-8}
 we announce the results of our computation of \mnn s
 for links with 6,7,8 crossings.
 The details will appear in the subsequent article.
  In Section \ref{s:irrat}
 we generalize the theorem \label{f:ineqq} for the case of Morse forms.
Section \ref{s:conj} contains some conjectures.
\subsection*{Acknowledgements}
The first author was supported by JST SPRING, Japan Grant Number JPMJSP2180.
The second author thanks the Nantes University, the DefiMaths program, and ``Mission Invite'' for the support and warm hospitality.
The second author was supported by JSPS KAKENHI Grant Numbers 24K06707, 20K03578.
The third author was supported by the  JSPS International Fellowships for Research in Japan (Short-term), FY2022 Research Abroad and Invitation Program for International Collaboration, and the WRH program of Tokyo Institute of Science. He thanks Tokyo Institute of Science for the support and warm hospitality.
 \newpage
--------------------

%@@@@@@@@@@@@@@@@@@@@@@@@@@@@@@@@@@@@
%@@@@@@@@@@@@@@@@@@@@@@@@@@@@@@@@@@@@@

Let $M$ be a compact manifold.
For a class $\xi\in H^1(M)$
the Morse-Novikov number
$\MM\NN(M,\xi)$
is defined as the minimal number of critical points of a
Morse map $f:M\to S^1$, such that the homotopy class of $f$
equals $\xi\in H^1(M)\approx [M, S^1]$ and
$f|\pr M$ is a fibration.
\mntt~ provides algebraic lower bounds for the \mnn.
These lower bounds, derived
 from Novikov homology, are in general
hard to compute
and can fail to be exact.

Consider the case when $M$ is the exterior $E(L)$
of a link $L$ in $S^3$.
Choose an orientation of $L$. For the cohomology class
$\Omega\in  H^1(E(L))$
induced by the orientation, the condition
$\MM\NN(E(L),\Omega)=0$
is equivalent to the condition of $L$ being fibered.
The question of determining of
$\MM\NN(E(L),\Omega)$
in the non-fibered case was addressed in several
articles.
%, see
%\cite{Goda}, \cite{PRW}, \cite{Manjarrez},
%\cite{Baker}.

In a series of papers
H. Goda gave a computation of Morse-Novikov numbers
$\MM\NN(E(L),\Omega)$
for  prime knots with at most 10 crossings
\cite{Goda}, \cite{Goda2006bis} and
links with at most  9 crossings \cite{Goda2006}
using the theory of sutured maniflds.

In the paper \cite{PRW} of
A. Pajitnov, L. Rudolph and C. Weber
the relations of Novikov homology and Alexander polynomial of links
was established. In the paper \cite{P-t} the third author proved that
\bq\label{f:citing-tunn}
\MM\NN(E(L),\Omega) \leq t(L)
\end{equation}
where $t(L)$ is the tunnel number of the link $L$.
%is less than or equal to twice the tunnel number of $L$.
The additivity of the number $\MM\NN(E(K),\Omega)$
\wrt~ connected sum of knots was established in the article
F. Manjarrez-Gutierrez \cite{Manjarrez} for the case of small knots,
and in the article of
K. Baker \cite{Baker} in general case.

%In the paper \cite{BakerManjarrezGutierrez} K. Baker and F. Manjarrez-Gutierrez
%introduced  an invariant $h(M,\g,\xi)$ for a sutured 3-manifold $(M,\g)$.

In the present paper we develop
algebraic and geometric tools, that allow us
to compute the number
$\MM\NN(E(L),\xi)$
for arbitrary classes $\xi$ in many cases.

%@@@@@@@@@@@@@@@@@@@@@@@@@@@@@@@@@@@@
%@@@@@@@@@@@@@@@@@@@@@@@@@@@@@@@@@@@@@

For the case of 3-manifolds the invariant $\MM\NN(M,\xi)$
is presumably equal to the invariant $h(M,\g,\xi)$ introduced
by K. Baker and F. Majarrez-Gutierrez
in \cite{BakerManjarrezGutierrez}, for some special choice of the suture $\g$.

%@@@@@@@@@@@@@@@@@@@@@@@@@
%@@@@@@@@@@@@@@@

%% file: auxxx.tex
\renewcommand{\a}{\alpha}
\renewcommand{\b}{\beta}
\newcommand{\g}{\gamma}

\newcommand{\e}{\epsilon}

\renewcommand{\t}{\theta}
\renewcommand{\l}{\lambda}

\newcommand{\m}{\mu}

\renewcommand{\r}{\rho}

\newcommand{\s}{\sigma}

\renewcommand{\o}{\omega}

\newcommand{\G}{\Gamma}
\newcommand{\D}{\Delta}

\renewcommand{\L}{\Lambda}

\renewcommand{\O}{\Omega}

\renewcommand{\AA}{{\mathcal A}}

\newcommand{\MM}{{\mathcal M}}
\newcommand{\NN}{{\mathcal N}}

\newcommand{\cc}{{\mathbb{C}}}

\newcommand{\nn}{{\mathbb{N}}}

\newcommand{\qq}{{\mathbb{Q}}}
\newcommand{\rr}{{\mathbb{R}}}

\newcommand{\zz}{{\mathbb{Z}}}

\newcommand{\AAAA}{{\mathscr{A}}}
\newcommand{\BBBB}{{\mathscr{B}} }
\newcommand{\CCCC}{{\mathscr{C}} }

\newcommand{\EEEE}{{\mathscr{E}} }

\newcommand{\OOOO}{{\mathscr{O}}}
\newcommand{\PPPP}{{\mathscr{P}}}

\newcommand{\UUUU}{{\mathscr{U}}}

\newcommand{\kkrest}{\begin{picture}(14,14)
\put(00,04){\line(1,0){14}}
\put(00,02){\line(1,0){14}}
\put(06,-4){\line(0,1){14}}
\put(08,-4){\line(0,1){14}}
\end{picture}     }

\newcommand{\krest}{~\kkrest~}

\newcommand{\Ker}{\text{\rm Ker }}

\newcommand{\rk}{\text{\rm rk}}
\renewcommand{\Im}{\text{\rm Im }}
\newcommand{\supp}{\text{\rm supp }}
\newcommand{\Int}{\text{\rm Int }}
\newcommand{\grad}{\text{\rm grad}}

\newcommand{\Id}{\text{\rm Id}}

\newcommand{\sgn}{\text{\rm sgn}}

\newcommand{\bere}{\begin{rema}}
\newcommand{\bede}{\begin{defi}}

\renewcommand{\beth}{\begin{theo}}
\newcommand{\bele}{\begin{lemm}}
\newcommand{\bepr}{\begin{prop}}
\newcommand{\beeq}{\begin{equation}}
\newcommand{\bega}{\begin{gather}}
\newcommand{\begaa}{\begin{gather*}}
\newcommand{\been}{\begin{enumerate}}

\newcommand{\bedee}{\begin{defii}}
\newcommand{\bethh}{\begin{theoo}}
\newcommand{\belee}{\begin{lemmm}}
\newcommand{\beprr}{\begin{propp}}

\newcommand{\beco}{\begin{coro}}

\newcommand{\beal}{\begin{aligned}}

\newcommand{\enre}{\end{rema}}

\newcommand{\enco}{\end{coro}}
\newcommand{\enpr}{\end{prop}}
\newcommand{\enth}{\end{theo}}
\newcommand{\enle}{\end{lemm}}
\newcommand{\enen}{\end{enumerate}}
\newcommand{\enga}{\end{gather}}
\newcommand{\engaa}{\end{gather*}}
\newcommand{\eneq}{\end{equation}}
\newcommand{\enal}{\end{aligned}}

\newcommand{\bq}{\begin{equation}}
\newcommand{\bqq}{\begin{equation*}}

\renewcommand{\leq}{\leqslant}
\renewcommand{\geq}{\geqslant}

\newcommand{\wi}{\widetilde}

\newcommand{\ove}{\overline}

\newcommand{\wh}{\widehat}

\newcommand{\sm}{\setminus}
\newcommand{\ems}{\varnothing}
\newcommand{\sbs}{\subset}

\newcommand{\tens}[1]{\underset{#1}{\otimes}}

\newcommand{\Lxi}{{\wh \L}_\xi}

\newcommand{\Geta}{{\wh \G}_\eta}

\newcommand{\Pxii}{{\wh \PPPP}_{\xi_i}}
\newcommand{\Axi}{{\wh \AAAA}_\xi}
\newcommand{\Bxi}{{\wh \BBBB}_\xi}

\newcommand{\Peta}{{\wh \PPPP}_{\eta}}

\newcommand{\ifff}{if and only if}

\newcommand{\sut}{~such~that~}

\newcommand{\wrt}{with respect to}
\newcommand{\ho}{homomorphism}

\newcommand{\ma}{manifold}
\newcommand{\nei}{neighbourhood}

\newcommand{\noconf}{~no~confusion~is~possible}

\newcommand{\Prf}{{\it Proof.\quad}}

\newcommand{\smo}{C^{\infty}}

\newcommand{\chart}{\Phi_p:U_p\to B^n(0,r_p)}
\newcommand{\atlas}{\{\Phi_p:U_p\to B^n(0,r_p)\}_{p\in S(f)}}

\newcommand{\pr}{\partial}

\newcommand{\prpr}[2]{\frac {\pr {#1}}{\pr {#2}}}

\newcommand{\qs}{\hfill\square}

\newcommand{\lLL}{\wh{\wh \L}}